\documentclass[final]{siamltex}
\usepackage{amsfonts}
\usepackage{amsmath}
\usepackage{amssymb}
\usepackage{mathtools}
\usepackage{array}
\usepackage{color}
\usepackage{eucal}
\usepackage{graphicx}
\usepackage{mathrsfs}
\usepackage{tikz}
\usepackage{rotating}
\usepackage[notref,notcite]{showkeys}

\newcommand{\curls}{\mathop{{\nabla\times}}}
\newcommand{\divs}{\mathop{\nabla\cdot}}
\newcommand{\grads}{\mathop{\nabla}}
\numberwithin{equation}{section}

\newcommand{\tr}{\mathrm{tr}}
\newcommand{\kerl}{\mathrm{Ker}}

\newcommand{\VT}{\mathcal{V}}
\newcommand{\EG}{\mathcal{E}}
\newcommand{\FC}{\mathcal{F}}
\newcommand{\vt}{\mathrm{v}}
\newcommand{\eg}{\mathrm{e}}
\newcommand{\fc}{\mathrm{f}}

\newcommand{\dK}{{\partial K}}
\newcommand{\HH}{{{C}^\infty}}

\newcommand{\oo}{\overset{\circ}}

\newcommand{\bld}[1]{\boldsymbol{#1}}
\newcommand{\bR}{\mathbb R}
\newcommand{\pol}{\mathcal{P}}
\newcommand{\bpol}{\bld{\mathcal{P}}}
\newcommand{\qol}{\mathcal{Q}}
\newcommand{\bqol}{\bld{\mathcal{Q}}}

\definecolor{red}{rgb}{0,0,0} 
\newcommand\gfu[1]{\textcolor{red}{#1}}
\definecolor{blue}{rgb}{0,0,0} 
\newcommand\jbc[1]{\textcolor{blue}{#1}}
\newcommand{\satri}{\mathrm{S_{1,k}^{
{\begin{tikzpicture}
\draw (0.1,0.1) -- (.1,.25)--(0.25,.1) -- (0.1,0.1);
\end{tikzpicture}}
}}}
\newcommand{\sbtri}{\mathrm{S_{2,k}^{
{\begin{tikzpicture}
\draw (0.1,0.1) -- (.1,.25)--(0.25,.1) -- (0.1,0.1);
\end{tikzpicture}}
}}}

\newcommand{\satrip}{\mathrm{S_{1,k+1}^{
{\begin{tikzpicture}
\draw (0.1,0.1) -- (.1,.25)--(0.25,.1) -- (0.1,0.1);
\end{tikzpicture}}
}}}

\newcommand{\sasqr}{\mathrm{S_{1,k}^{
{\begin{tikzpicture}
\draw (0.1,0.1) -- (.1,.25)-- (0.25,0.25) -- (0.25,.1) -- (0.1,0.1);
\end{tikzpicture}}
}}}
\newcommand{\sbsqr}{\mathrm{S_{2,k}^{
{\begin{tikzpicture}
\draw (0.1,0.1) -- (.1,.25)-- (0.25,0.25) -- (0.25,.1) -- (0.1,0.1);
\end{tikzpicture}}
}}}
\newcommand{\scsqr}{\mathrm{S_{3,k}^{
{\begin{tikzpicture}
\draw (0.1,0.1) -- (.1,.25)-- (0.25,0.25) -- (0.25,.1) -- (0.1,0.1);
\end{tikzpicture}}
}}}
\newcommand{\sdsqr}{\mathrm{S_{4,k}^{
{\begin{tikzpicture}
\draw (0.1,0.1) -- (.1,.25)-- (0.25,0.25) -- (0.25,.1) -- (0.1,0.1);
\end{tikzpicture}}
}}}
\newcommand{\sisqr}{\mathrm{S_{i,k}^{
{\begin{tikzpicture}
\draw (0.1,0.1) -- (.1,.25)-- (0.25,0.25) -- (0.25,.1) -- (0.1,0.1);
\end{tikzpicture}}
}}}

\newcommand{\sasqrp}{\mathrm{S_{1,k+1}^{
{\begin{tikzpicture}
\draw (0.1,0.1) -- (.1,.25)-- (0.25,0.25) -- (0.25,.1) -- (0.1,0.1);
\end{tikzpicture}}
}}}

\newcommand{\satet}{
\mathrm{S_{1,k}^{
 {
\begin{tikzpicture}[scale=0.4]
\draw  (0,0) -- (0,0.3);
\draw  (0,0) -- (0.2996,-0.0157);
\draw  (0,0) -- (-.2194,-0.2046);
\draw (0,0.3) -- (-.2194,-0.2046)--(0.2996,-0.0157)--(0,0.3);
\end{tikzpicture}}
}}}

\newcommand{\sbtet}{
\mathrm{S_{2,k}^{
 {
\begin{tikzpicture}[scale=0.4]
\draw  (0,0) -- (0,0.3);
\draw  (0,0) -- (0.2996,-0.0157);
\draw  (0,0) -- (-.2194,-0.2046);
\draw (0,0.3) -- (-.2194,-0.2046)--(0.2996,-0.0157)--(0,0.3);
\end{tikzpicture}}
}}}

\newcommand{\sacub}{\mathrm{S_{1,k}^{
\begin{tikzpicture}[scale=0.4]
\draw  (0,0) -- (0,0.3);
\draw  (0,0) -- (0.2996,-0.0157);
\draw  (0,0) -- (-.2194,-0.2046);
\draw (-.2194,-0.2046) -- (.0802,-0.2046) -- (0.2996,-0.0157);
\draw  (0,0.3) -- (0.2996,0.2843);
\draw (0,0.3) -- (-.2194,0.0954);
\draw (-.2194,0.0954) -- (.0802,0.0954) -- (0.2996,0.2843);
\draw (-.2194,0.0954) -- (-.2194,-0.2046);
\draw (.0802,0.0954) -- (.0802,-0.2046);
\draw (0.2996,-0.0157) -- (0.2996,0.2843);
\end{tikzpicture}
}}}

\newcommand{\sbcub}{\mathrm{S_{2,k}^{
\begin{tikzpicture}[scale=0.4]
\draw  (0,0) -- (0,0.3);
\draw  (0,0) -- (0.2996,-0.0157);
\draw  (0,0) -- (-.2194,-0.2046);
\draw (-.2194,-0.2046) -- (.0802,-0.2046) -- (0.2996,-0.0157);
\draw  (0,0.3) -- (0.2996,0.2843);
\draw (0,0.3) -- (-.2194,0.0954);
\draw (-.2194,0.0954) -- (.0802,0.0954) -- (0.2996,0.2843);
\draw (-.2194,0.0954) -- (-.2194,-0.2046);
\draw (.0802,0.0954) -- (.0802,-0.2046);
\draw (0.2996,-0.0157) -- (0.2996,0.2843);
\end{tikzpicture}
}}}
\newcommand{\sccub}{\mathrm{S_{3,k}^{
\begin{tikzpicture}[scale=0.4]
\draw  (0,0) -- (0,0.3);
\draw  (0,0) -- (0.2996,-0.0157);
\draw  (0,0) -- (-.2194,-0.2046);
\draw (-.2194,-0.2046) -- (.0802,-0.2046) -- (0.2996,-0.0157);
\draw  (0,0.3) -- (0.2996,0.2843);
\draw (0,0.3) -- (-.2194,0.0954);
\draw (-.2194,0.0954) -- (.0802,0.0954) -- (0.2996,0.2843);
\draw (-.2194,0.0954) -- (-.2194,-0.2046);
\draw (.0802,0.0954) -- (.0802,-0.2046);
\draw (0.2996,-0.0157) -- (0.2996,0.2843);
\end{tikzpicture}
}}}
\newcommand{\sdcub}{\mathrm{S_{4,k}^{
\begin{tikzpicture}[scale=0.4]
\draw  (0,0) -- (0,0.3);
\draw  (0,0) -- (0.2996,-0.0157);
\draw  (0,0) -- (-.2194,-0.2046);
\draw (-.2194,-0.2046) -- (.0802,-0.2046) -- (0.2996,-0.0157);
\draw  (0,0.3) -- (0.2996,0.2843);
\draw (0,0.3) -- (-.2194,0.0954);
\draw (-.2194,0.0954) -- (.0802,0.0954) -- (0.2996,0.2843);
\draw (-.2194,0.0954) -- (-.2194,-0.2046);
\draw (.0802,0.0954) -- (.0802,-0.2046);
\draw (0.2996,-0.0157) -- (0.2996,0.2843);
\end{tikzpicture}
}}}
\newcommand{\sicub}{\mathrm{S_{i,k}^{
\begin{tikzpicture}[scale=0.4]
\draw  (0,0) -- (0,0.3);
\draw  (0,0) -- (0.2996,-0.0157);
\draw  (0,0) -- (-.2194,-0.2046);
\draw (-.2194,-0.2046) -- (.0802,-0.2046) -- (0.2996,-0.0157);
\draw  (0,0.3) -- (0.2996,0.2843);
\draw (0,0.3) -- (-.2194,0.0954);
\draw (-.2194,0.0954) -- (.0802,0.0954) -- (0.2996,0.2843);
\draw (-.2194,0.0954) -- (-.2194,-0.2046);
\draw (.0802,0.0954) -- (.0802,-0.2046);
\draw (0.2996,-0.0157) -- (0.2996,0.2843);
\end{tikzpicture}
}}}

\newcommand{\sapri}{\mathrm{S_{1,k}^{
\begin{tikzpicture}[scale=0.4]
\draw  (0,0) -- (0,0.3);
\draw  (0,0) -- (0.2996,-0.0157);
\draw  (0,0) -- (-.2194,-0.2046);
\draw (-.2194,-0.2046) -- (0.2996,-0.0157);
\draw  (0,0.3) -- (0.2996,0.2843);
\draw (0,0.3) -- (-.2194,0.0954);
\draw (-.2194,0.0954)-- (0.2996,0.2843);
\draw (-.2194,0.0954) -- (-.2194,-0.2046);
\draw (0.2996,-0.0157) -- (0.2996,0.2843);
\end{tikzpicture}}
}}
\newcommand{\sbpri}{\mathrm{S_{2,k}^{
\begin{tikzpicture}[scale=0.4]
\draw  (0,0) -- (0,0.3);
\draw  (0,0) -- (0.2996,-0.0157);
\draw  (0,0) -- (-.2194,-0.2046);
\draw (-.2194,-0.2046) -- (0.2996,-0.0157);
\draw  (0,0.3) -- (0.2996,0.2843);
\draw (0,0.3) -- (-.2194,0.0954);
\draw (-.2194,0.0954)-- (0.2996,0.2843);
\draw (-.2194,0.0954) -- (-.2194,-0.2046);
\draw (0.2996,-0.0157) -- (0.2996,0.2843);
\end{tikzpicture}}
}}
\newcommand{\scpri}{\mathrm{S_{3,k}^{
\begin{tikzpicture}[scale=0.4]
\draw  (0,0) -- (0,0.3);
\draw  (0,0) -- (0.2996,-0.0157);
\draw  (0,0) -- (-.2194,-0.2046);
\draw (-.2194,-0.2046) -- (0.2996,-0.0157);
\draw  (0,0.3) -- (0.2996,0.2843);
\draw (0,0.3) -- (-.2194,0.0954);
\draw (-.2194,0.0954)-- (0.2996,0.2843);
\draw (-.2194,0.0954) -- (-.2194,-0.2046);
\draw (0.2996,-0.0157) -- (0.2996,0.2843);
\end{tikzpicture}}
}}
\newcommand{\sdpri}{\mathrm{S_{4,k}^{
\begin{tikzpicture}[scale=0.4]
\draw  (0,0) -- (0,0.3);
\draw  (0,0) -- (0.2996,-0.0157);
\draw  (0,0) -- (-.2194,-0.2046);
\draw (-.2194,-0.2046) -- (0.2996,-0.0157);
\draw  (0,0.3) -- (0.2996,0.2843);
\draw (0,0.3) -- (-.2194,0.0954);
\draw (-.2194,0.0954)-- (0.2996,0.2843);
\draw (-.2194,0.0954) -- (-.2194,-0.2046);
\draw (0.2996,-0.0157) -- (0.2996,0.2843);
\end{tikzpicture}}
}}

\newcommand{\sapyr}{\mathrm{S_{1,k}^{
\begin{tikzpicture}[scale=0.4]
\draw  (0,0) -- (0,0.3);
\draw  (0,0) -- (0.2996,-0.0157);
\draw  (0,0) -- (-.2194,-0.2046);
\draw (-.2194,-0.2046) -- (.0802,-0.2046) -- (0.2996,-0.0157);
\draw (0,0.3) -- (-.2194,-0.2046);
\draw (0,0.3) -- (.0802,-0.2046);
\draw (0,0.3) -- (0.2996,-0.0157);
\end{tikzpicture}}
}}
\newcommand{\sbpyr}{\mathrm{S_{2,k}^{
\begin{tikzpicture}[scale=0.4]
\draw  (0,0) -- (0,0.3);
\draw  (0,0) -- (0.2996,-0.0157);
\draw  (0,0) -- (-.2194,-0.2046);
\draw (-.2194,-0.2046) -- (.0802,-0.2046) -- (0.2996,-0.0157);
\draw (0,0.3) -- (-.2194,-0.2046);
\draw (0,0.3) -- (.0802,-0.2046);
\draw (0,0.3) -- (0.2996,-0.0157);
\end{tikzpicture}}
}}
\newcommand{\scpyr}{\mathrm{S_{3,k}^{
\begin{tikzpicture}[scale=0.4]
\draw  (0,0) -- (0,0.3);
\draw  (0,0) -- (0.2996,-0.0157);
\draw  (0,0) -- (-.2194,-0.2046);
\draw (-.2194,-0.2046) -- (.0802,-0.2046) -- (0.2996,-0.0157);
\draw (0,0.3) -- (-.2194,-0.2046);
\draw (0,0.3) -- (.0802,-0.2046);
\draw (0,0.3) -- (0.2996,-0.0157);
\end{tikzpicture}}
}}
\newcommand{\sdpyr}{\mathrm{S_{4,k}^{
\begin{tikzpicture}[scale=0.4]
\draw  (0,0) -- (0,0.3);
\draw  (0,0) -- (0.2996,-0.0157);
\draw  (0,0) -- (-.2194,-0.2046);
\draw (-.2194,-0.2046) -- (.0802,-0.2046) -- (0.2996,-0.0157);
\draw (0,0.3) -- (-.2194,-0.2046);
\draw (0,0.3) -- (.0802,-0.2046);
\draw (0,0.3) -- (0.2996,-0.0157);
\end{tikzpicture}}
}}
\newcommand{\sipyr}{\mathrm{S_{i,k}^{
\begin{tikzpicture}[scale=0.4]
\draw  (0,0) -- (0,0.3);
\draw  (0,0) -- (0.2996,-0.0157);
\draw  (0,0) -- (-.2194,-0.2046);
\draw (-.2194,-0.2046) -- (.0802,-0.2046) -- (0.2996,-0.0157);
\draw (0,0.3) -- (-.2194,-0.2046);
\draw (0,0.3) -- (.0802,-0.2046);
\draw (0,0.3) -- (0.2996,-0.0157);
\end{tikzpicture}}
}}

\title{A systematic construction of finite element commuting exact sequences}

\author{
Bernardo Cockburn
        \thanks{School of Mathematics, University of Minnesota, Vincent Hall,
                Minneapolis, MN 55455, USA, email: {\tt cockburn@math.umn.edu}.
                Supported in part by the National Science Foundation
                (Grant DMS-1522657) and by the University of Minnesota
                Supercomputing Institute.}
\and
Guosheng Fu
        \thanks{School of Mathematics, University of Minnesota, Vincent Hall,
                Minneapolis, MN 55455, USA, email: {\tt fuxxx165@math.umn.edu}.}
}

\date{\today}

\begin{document}


%

\maketitle

\begin{abstract}  
We present a systematic construction of finite element exact sequences with a commuting diagram 
for the de Rham complex in one-, two- and three-space dimensions. 
We apply the  construction in two-space dimensions
to rediscover two families of exact sequences for triangles and three for squares, and to 
uncover one new family of exact sequence for squares and 
two new families of exact sequences for general polygonal elements. 
We apply the construction in three-space dimensions to rediscover
two families of exact sequences for tetrahedra, three for cubes, and 
one for prisms; and to uncover four new families of exact sequences for pyramids,
three for prisms, and one for cubes. 
\end{abstract}



\thispagestyle{plain} \markboth{B. Cockburn and G. Fu}{
A systematic construction of commuting exact sequences}

\

\centerline{{\bf Version of \today}}

\



\section{Introduction}
\label{sec:introduction}
We give a systematic 
construction of finite element exact sequences 
with a commuting diagram property
for the de Rham complex in one-, two- and three-dimensional domains. 
The construction of these commuting exact sequences
for the three-dimensional case relies on the construction of commuting exact sequences for the two dimensional case, which in turn relies on the construction of commuting exact sequences in one-space dimension. For each dimension, the construction is carried out on a single element $K$ 
in such  way that proper continuity properties hold which guarantee that the commutativity of the diagram holds when $K$ is replaced by the triangulation of the domain $\Omega$, $\Omega_h:=\{K\}$. In two dimensions, the elements $K$ can be polygons of arbitrary shape, and in three-dimensions, polyhedra of arbitrary shape.
 
Roughly speaking, by a commuting exact sequence on a polyhedral domain $\Omega$ in $\bR^3$, we mean that, the mappings 
in the diagram
\begin{alignat*}{2}
\begin{tabular}{c c c c c c c}
$ \HH(\overline{\Omega})\;\;$&
$\overset{\grads}{\longrightarrow} $&
$\bld\HH(\overline{\Omega})\;\;$&
$\overset{\curls}{\longrightarrow} $&
$\bld\HH(\overline{\Omega})\;\;$&
$\overset{\divs}{\longrightarrow} $&
$\HH(\overline{\Omega})\;\;$\\
$\downarrow {\scriptstyle \Pi_H}$ & &
$\downarrow{\scriptstyle \Pi_E}$
& &$\downarrow{\scriptstyle \Pi_V}$
& &
$\downarrow{\scriptstyle \Pi_W}$
\\
$H(\Omega_h)\;\;$&
$\overset{\grads}{\longrightarrow}$ &
$E(\Omega_h)\;\;$&
$\overset{\curls}{\longrightarrow}$ &
$V(\Omega_h)\;\;$&
$\overset{\divs}{\longrightarrow}$ &
$W(\Omega_h)\;\;$
\end{tabular}
\end{alignat*}
commute in the sense that 
\begin{alignat*}{2}
 \grads \Pi_H  u = &\;\Pi_E \grads u&&\;\;\;\;\forall u\in \HH(\overline{\Omega}),\\
 \curls \Pi_E  u = &\; \Pi_V \curls u&&\;\;\;\;\forall u\in \bld \HH(\overline{\Omega}),\\
 \divs \Pi_V  u = &\; \Pi_W \divs u&&\;\;\;\;\forall u\in \bld \HH(\overline{\Omega}).
\end{alignat*}
The importance of commuting exact  sequences, for the devising of stable 
finite element methods, has been amply discussed in 
\cite{ArnoldFalkWinther06, ArnoldFalkWinther10} and the references therein. 

Most of the previous work on the construction of commuting exact sequences
(in three-space dimensions) focuses on the explicit construction of shape functions
on one of {four} particular {\it reference} polyhedra, namely, the 
reference tetrahedron, hexahedron (cube), prism, and pyramid. 
See \cite{Nedelec80,Nedelec86} for sequences on the reference tetrahedron and
reference hexahedron,  \cite{ArnoldAwanou14,CockburnQiu14} on the reference hexahedron,  
\cite{Zaglmayr06} on the { reference} prism, and \cite{NigamPhillips12,NigamPhillips12b} on the { reference}
pyramid.  All of these spaces in these sequences are spanned by polynomial shape functions, 
except those in  \cite{NigamPhillips12,NigamPhillips12b} which also contain rational shape functions.

Currently, there are two ways of constructing commuting exact sequences on general polyhedral meshes.
The first is provided by the {\em Virtual Element Methods} (VEM)\jbc{; see the 2013 paper  \cite{VeigaBrezziCangianiManziniMariniRusso13} and the 2015 paper
\cite{VeigaBrezziMariniRusso15}.} 
These methods define the basis functions of their local spaces on each polyhedral element in terms of solutions to 
certain partial differential equations.
The explicit form of these basis functions, usually not computable, is {\em not} needed by the methods, but 
a set of {\em unisolvent} degrees of freedom which can be used  to {\em exactly} compute integrals related to 
the polynomial parts of the basis functions needs be constructed. In our construction, we avoid basis functions defined as solutions to certain partial differential equations.

The second way is provided by the so-called {\em Finite Element System} (FES) developed \jbc{in 2011} in \cite[Section 5]
{ChristiansenMuntheKaasOwren11}; see also \jbc{the recent papers} \cite{ChristiansenRapetti16,ChristiansenGillette15}. Therein, the notion 
of a {\em compatible} FES was introduced, see \cite[Definition 5.12]{ChristiansenMuntheKaasOwren11},  which was 
then proven to \gfu{be equivalent to} the existence of a commuting diagram, 
see \cite[Proposition 5.44]{ChristiansenMuntheKaasOwren11}.
The construction of a FES with a commuting diagram was thus reduced to the construction of a compatible FES.
\gfu{In \cite[Example 5.29]{ChristiansenMuntheKaasOwren11},  the authors obtained 
compatible FESs on a general ($n$-dimensional) polytope mesh via {\em element agglomeration} from (available) compatible FESs defined on a refined (simplicial) mesh. 
However, the resulting compatible FES on the polytope mesh is not particularly interesting
since it either provides identical spaces on the refined (simplicial) mesh or requires finding a subsystem of  {\em locally harmonic forms} which seems to be a very hard task.}

Furthermore, in \cite[Corollary 3.2]{ChristiansenGillette15}, a criterion for finding a compatible FES, containing certain prescribed  functions, with the {\em smallest} possible dimension was presented given thus rise to the concept of a {\em minimal compatible FES} (mcFES). 
\gfu{\jbc{The minimality was proven} for three finite element systems, namely, the  {\em trimmed polynomial differential forms} \cite{ArnoldFalkWinther06} on a simplicial mesh (in three-dimensions \jbc{this} is the exact sequence due to Ned\'el\'ec \cite{Nedelec80}), the serendipity elements \cite{ArnoldAwanou14} on a cubic mesh,
and the TNT elements \cite{CockburnQiu14} on a cubic mesh; these seem to be the only mcFESs available in the literature.} 
Although a simple dimension count equation in \cite[Corollary 3.2]{ChristiansenGillette15} 
can be used to check whether a compatible FES is minimal or not, a  practical construction of an mcFES was not provided  in \cite{ChristiansenGillette15}. 

\gfu{Our work can be naturally considered as a practical way to construct an mcFES  in  one-, two- and three-space 
dimensions.} 
We develop a systematic construction of commuting exact sequences on a general polytope
in \gfu{one-, two- and three-space 
dimensions},  and then (in three-space dimensions) apply the construction to each of the above-mentioned four {reference} polyhedra to discover and rediscover concrete examples of commuting exact sequences. In particular, we show that all the results in 
\cite{Nedelec80,Nedelec86,ArnoldAwanou14,CockburnQiu14,Zaglmayr06}
on a tetrahedron, cube, and prism fit nicely within our construction. Moreover, a 
significant dimension reduction is obtained on the pyramidal commuting exact sequence in comparison with
the exact sequence proposed in \cite{NigamPhillips12,NigamPhillips12b}.
For a general polyhedron, our construction of of high-order commuting exact sequences  is significantly more difficult than that of the cases already mentioned. It will be carried out elsewhere.

Let us now roughly describe the main steps of the systematic construction in three-dimensions. We proceed as follows.
First, we pick $K$ to be a specific polyhedron. Then, we take the set of commuting exact sequences on the faces of $K$,
\begin{align*}
S_2(\partial K):=\{ S_2(\fc):\;\fc\in\FC(K)\}. 
\end{align*}
Each of the commuting exact sequences 

\

\begin{tabular}{c c c  c c c c c}
$\phantom{oooooooooooooo}S_2(\fc): $&$H_2(\fc)$&
$\overset{\grads}{\longrightarrow}$ &
$E_2(\fc)$
&
$\overset{\curls}{\longrightarrow}$ &
$W_2(\fc)$,
\end{tabular}

\

\noindent was {\em previously} obtained by applying the systematic construction in the two dimensional case. 
Moreover, we require a 
{\em compatibility} condition on the edges, namely, that, if the faces $\fc_1$ and $\fc_2$ share the edge $\eg$, the trace on the edge $\eg$
of $H_2(\fc_1)$ must coincide with that of $H_2(\fc_2)$, and that the tangential trace on the edge $\eg$
of $E_2(\fc_1)$ must coincide with that of $E_2(\fc_2)$.
Next, we consider a candidate for the commuting exact  sequence we are seeking of the form

\

\begin{tabular}{c c c  c c c c c}
$\phantom{oooooooo}S_3^g(K): $&$H_3^g(K)$&
$\overset{\grads}{\longrightarrow}$ &
$E_3^g(K)$
&
$\overset{\curls}{\longrightarrow}$ &
$V_3^g(K)$
&
$\overset{\divs}{\longrightarrow}$ &
$W_3^g(K)$,
\end{tabular}

\

\noindent We require $S^g_3(K)$ to be exact; to have, for each face $\fc$ of $K$, the traces of $H_3^g(K)$ 
on the face $\fc$  in the space $H_2(\fc)$, the traces of $E_3^g(K)$ on
the face $\fc$  in the space $E_2(\fc)$, and the traces of $V_3^g(K)$ on
the face $\fc$  in the space $W_2(\fc)$; and to have the constant functions in the space  $W_3^g(K)$.

Finally, we {\em characterize} the spaces $\delta H_3(K)$ and $\delta E_3(K)$ as the spaces with smallest dimension 
such that

\begin{align*}
\begin{tabular}{c c c c c c c c}
&$H_3^g(K)$&
$\overset{\grads}{\longrightarrow}$ &
$E_3^g(K)$
&
$\overset{\curls}{\longrightarrow}$ &
$V_3^g(K)$
&
$\overset{\divs}{\longrightarrow}$ &
$W_3^g(K)$
\\
&$\oplus \delta H_3(K)$&
&
$\oplus\grads\delta H_3(K)$
&&
$\oplus\curls\delta E_3(K)$
&
&
\\
&&
&
$\oplus \delta E_3(K)$
&
&&
&
\end{tabular}
\end{align*}
is a commuting exact sequence for which, for each face $\fc$ of $K$, the traces of $H_3^g(K)\oplus \delta H_3(K)$ 
on the face $\fc$ constitute the space $H_2(\fc)$, the tangential traces of $E_3^g(K)\oplus \grads\delta H_3(K)\oplus\delta E(K)$ on
the face $\fc$ constitute the space $E_2(\fc)$, and the normal traces of $V_3^g(K)\oplus \curls\delta E_3(K)$ on
the face $\fc$ constitute the space $W_2(\fc)$.

This completes the rough description of the systematic construction. Note that the spaces $\delta H_3(K)$ and $\delta E_3(K)$
are {\em not} necessarily unique even though they do have their dimension is uniquely determined.  
Note also that the construction of the space $\delta E_3(K)$, which is the {\em most difficult} 
part of the general construction, is essentially a particular case of the construction of $M$-decompositions for mixed methods, see \cite{CockburnFuSayas16, CockburnFuM2D,CockburnFuM3D}.  Thus, the present  construction can be considered to be an extension of the $M$-decomposition approach to 
the setting of exact sequences for the de Rham complex.

\begin{table}[!ht]
 \caption{Commuting exact  sequences obtained as the application of the systematic construction. 
 The index $k$ is associated to a polynomial degree.  The symbol $\checkmark$ indicates that the sequence is new.
 Sequences in the same row are \jbc{such} that the spaces of traces are the same on faces of the same shape.}
\centering
\begin{tabular}{c c c c c }
\hline
\noalign{\smallskip} 
sequence &
${
\begin{tikzpicture}
\draw  (0,0) -- (0,0.3);
\draw  (0,0) -- (0.2996,-0.0157);
\draw  (0,0) -- (-.2194,-0.2046);
\draw (0,0.3) -- (-.2194,-0.2046)--(0.2996,-0.0157)--(0,0.3);
\end{tikzpicture}}$
&
 ${
\begin{tikzpicture}
\draw  (0,0) -- (0,0.3);
\draw  (0,0) -- (0.2996,-0.0157);
\draw  (0,0) -- (-.2194,-0.2046);
\draw (-.2194,-0.2046) -- (.0802,-0.2046) -- (0.2996,-0.0157);
\draw  (0,0.3) -- (0.2996,0.2843);
\draw (0,0.3) -- (-.2194,0.0954);
\draw (-.2194,0.0954) -- (.0802,0.0954) -- (0.2996,0.2843);
\draw (-.2194,0.0954) -- (-.2194,-0.2046);
\draw (.0802,0.0954) -- (.0802,-0.2046);
\draw (0.2996,-0.0157) -- (0.2996,0.2843);
\end{tikzpicture}}$ 
&
 ${
\begin{tikzpicture}
\draw  (0,0) -- (0,0.3);
\draw  (0,0) -- (0.2996,-0.0157);
\draw  (0,0) -- (-.2194,-0.2046);
\draw (-.2194,-0.2046) -- (0.2996,-0.0157);
\draw  (0,0.3) -- (0.2996,0.2843);
\draw (0,0.3) -- (-.2194,0.0954);
\draw (-.2194,0.0954)-- (0.2996,0.2843);
\draw (-.2194,0.0954) -- (-.2194,-0.2046);
\draw (0.2996,-0.0157) -- (0.2996,0.2843);
\end{tikzpicture}}$ 
&
  ${
\begin{tikzpicture}
\draw  (0,0) -- (0,0.3);
\draw  (0,0) -- (0.2996,-0.0157);
\draw  (0,0) -- (-.2194,-0.2046);
\draw (-.2194,-0.2046) -- (.0802,-0.2046) -- (0.2996,-0.0157);
\draw (0,0.3) -- (-.2194,-0.2046);
\draw (0,0.3) -- (.0802,-0.2046);
\draw (0,0.3) -- (0.2996,-0.0157);
\end{tikzpicture}}$ 
\\
\noalign{\smallskip}
\hline\hline
\noalign{\smallskip} 
${\mathrm{S_{1,k}}}$ 
& 
\cite{Nedelec86} & \cite{ArnoldAwanou14} & $\checkmark$& $\checkmark$
\\
$\mathrm{S_{2,k}}$ & 
\cite{Nedelec80} & $\checkmark$& $\checkmark$& $\checkmark$
\\
\noalign{\smallskip} 
\hline
\noalign{\smallskip} 
$\mathrm{S_{3,k}}$ & 
- & \cite{CockburnQiu14} & $\checkmark$& $\checkmark$
\\
$\mathrm{S_{4,k}}$ & 
- & \cite{Nedelec80}& \cite{FuentesKeithDemkowiczNagaraj15}& $\checkmark$
\\
\noalign{\smallskip} 
\hline 
\\
\end{tabular}
\label{table:sequences}
\end{table}

Let us now comment on the actual commuting exact  sequences we obtain here.
We apply the systematic construction to explicitly obtain {\it fourteen} families of 
commuting exact sequences:  on  
the reference tetrahedron (2 sequences), reference cube (4 sequences), 
reference prism (4 sequences), and reference pyramid (4 sequences).  
In Table \ref{table:sequences}, 
we indicate if these are known or new commuting exact  
sequences. Let us note that, on a reference prism, 
there is an additional family of commuting exact  sequence, namely, 
the one proposed in  \cite{Zaglmayr06}; 
it uses the $H(\mathrm{div})$ and $L^2$ spaces obtained in \cite{ChenDouglas89}. The sequence introduced in  \cite{FuentesKeithDemkowiczNagaraj15} is a slight modification of this one. 
On a reference pyramid, there are two additional 
family of exact sequences in \cite{NigamPhillips12,NigamPhillips12b}. 
Our spaces are significantly smaller.

Our {fourteen} commuting exact  sequences can be gathered into  {\it four} groups of sequences each of which is
displayed in a row in Table \ref{table:sequences}. These {\it four} group of sequences are \jbc{such}  
that the $H^1$-, $H(\mathrm{curl})$-, and $H(\mathrm{div})$- trace spaces on  similar faces are the {\em same}.
As a consequence, each of the {\it four} groups of sequences 
can be patched into in a {\it hybrid polyhedral mesh} $\Omega_h=\{K\}$ where the elements $K$ are suitably defined
by affine mappings of the four reference polyhedral elements. This way of regrouping the sequences is motivated by the work proposed in 
\cite{FuentesKeithDemkowiczNagaraj15},  where one group of sequences that can be used on a hybrid mesh
(with a more complicated  pyramidal sequence from \cite{NigamPhillips12}) 
was recently carefully studied to construct ``orientation embedded high-order shape functions''.

The rest of the paper is organized as follows. In Section \ref{sec:mainresults}, 
we 
present our main results on the systematic construction of commuting exact sequences
in one-, two- and three-space dimensions. 
Then in Section \ref{sec:applications}, we apply the systematic construction to 
explicitly obtain commuting exact sequences
on the reference interval in one-space dimension;
on the reference triangle, reference square, and 
on a general polygon in two-space dimensions;
and on the above-mentioned four reference polyhedra 
in three-space dimensions.
Section \ref{sec:proof} is devoted to the proofs of the
Section \ref{sec:proof2} is devoted to the proofs results in  
Section \ref{sec:applications}.
We end in Section \ref{sec:conclude} with some concluding remarks.

\section{A systematic construction of commuting exact sequences}
\label{sec:mainresults}
In this section, we introduce the notation used throughout the paper. We then define the concept of a 
{\em compatible exact sequence} in one-, two- and  three-space dimensions, which, in differential form language, is nothing but the {\em compatible} FES introduced in 
\cite[Definition 5.12]{ChristiansenMuntheKaasOwren11}.  Let us recall that, the theory of FESs introduced in \cite[Proposition 5.44]{ChristiansenMuntheKaasOwren11}, 
establishes the equivalence of a compatible exact sequence and a sequence admitting a commuting diagram. As a consequence, the construction of a sequence admitting a commuting diagram is reduced to the construction of a compatible exact sequence. We use this powerful result and devote ourselves to developing a systematic construction of compatible exact sequences in one-, two- and three-space dimensions. 

In Section \ref{sec:applications}, this approach is applied to obtain many compatible exact sequences 
with explicitly-defined shape functions for elements of various shapes.

\subsection{Notation}
Here we introduce the notation we use for the rest of the paper.

\subsection*{Geometry}
We denote by  $K\subset \bR^d$, a segment if $d=1$, a polygon if $d=2$, and a polyhedron if $d=3$. 
We denote its boundary by $\dK$, the set of its vertices by $\VT(K)$, the set of its edges (for $d=1,2$) by $\EG(K)$,
and the set of its faces (for $d=3$) by $\FC(K)$. 

\subsection*{Trace operators}
For a scalar-valued function $v$ on $K$ with sufficient regularity,
we denote by $\tr_H^\vt v:= v|_\vt$ the trace of $v$ on a vertex $\vt\in \VT(K)$, 
by $\tr_H^\eg v:= v|_\eg$ the trace of $v$ on an edge $\eg\in \EG(K)$, 
by $\tr_H^\fc v:= v|_\fc$ the trace of $v$ on a face $\fc\in \FC(K)$, and  
by $\tr_H v :=v|_\dK$ the trace of $v$ on the whole boundary $\dK$.

If $d\ge 2$, for a $d-$dimensional vector-valued function $v$ with sufficient regularity,
we denote by $\tr_E^\eg v:= \left.(v\cdot t_{\eg})\right|_\eg$, where $t_{\eg}$ is the unit vector in the direction of the edge $\eg$,  the tangential 
trace of $v$ on an edge $\eg\in \EG(K)$. We denote  by  
$\tr_E^\fc v:= \left.\left(n_\fc\times(v\times n_\fc)\right)\right|_\fc$, where $n_{\fc}$ is the unit outward normal to the face $\fc$,  the tangential 
trace of $v$ on a face $\fc\in \FC(K)$. Finally, we denote
by 
$\tr_E v:= \left\{\begin{tabular}{l l}
                       $\left.(v\cdot t_\dK)\right|_\dK$ & if $d = 2$,\\
                       $\left.\left(n_\dK\times(v\times n_\dK)\right)\right|_\dK$ & if $d = 3$,
                      \end{tabular}\right.
$ the tangential 
trace of $v$ on the whole boundary $\dK$. Here, $t_\dK$ on the edge $\eg$ is nothing but $t_\eg$.
Similarly,  $n_\dK$ on the face $\fc$ is nothing but $n_\fc$.

Finally, if $d= 3$, for a $d-$dimensional vector-valued  function $v$ with sufficient regularity,
we denote by 
$\tr_V^\fc v:= \left.\left(v\cdot n_\fc\right)\right|_\fc$ the normal 
trace of $v$ on a face $\fc\in \FC(K)$, and by 
$\tr_V v:= \left.\left(v\cdot n_\dK\right)\right|_\dK$ the normal 
trace of $v$ on the whole boundary $\dK$.


\subsection*{Differential operators}
If $d = 1$, we denote by $x_1$ the coordinate of a point on $K$, 
if $d = 2$, we denote by $(x_1, x_2)$  the coordinate of a point on $K$,
and if $d = 3$, we denote by $(x_1,x_2,x_3)$ the coordinate of a point on $K$. 
We define the gradient operator $\grads$, curl operator $\curls$ (for $d\ge 2$), 
and divergence operator $\divs$ (for $d=3$) on the element $K$
as follows.

For a scalar function $v$ on $K$ with sufficient regularity, we set
\begin{align*}
 \grads v:=\left\{
 \begin{tabular}{l l}
  $\partial_{x_1} v$& if $d=1$,\\
  $(\partial_{x_1} v, \partial_{x_2} v)^t$& if $d=2$,\\
  $(\partial_{x_1} v, \partial_{x_2} v, \partial_{x_3} v)^t$ & if $d=3$.
 \end{tabular}
 \right.
\end{align*}

For a $d$-dimensional vector function $v$ on $K$ 
with sufficient regularity ($d\ge 2$), we set
\begin{align*}
 \curls v:=\left\{
 \begin{tabular}{l l}
  $-\partial_{x_2} v_1+ \partial_{x_1} v_2$& if $d=2$,
  \vspace{.2cm}\\
  $\left(
  \begin{tabular}{c}
   $-\partial_{x_3} v_2+\partial_{x_2} v_3$\\
  $-\partial_{x_1} v_3+\partial_{x_3} v_1$\\
  $-\partial_{x_2} v_1+\partial_{x_1} v_2$ 
  \end{tabular} 
  \right)$ & if $d=3$,
 \end{tabular}
 \right.
\end{align*}
where $v=(v_1,v_2)^t$ if $d=2$ and $v=(v_1,v_2,v_3)^t$ if $d=3$.

For a $3$-dimensional vector function $v=(v_1,v_2,v_3)^t$ on $K$ 
with sufficient regularity, we set
\begin{align*}
 \divs v:= \partial_{x_1}v_1 +
 \partial_{x_2}v_2 +
 \partial_{x_3}v_3.
 \end{align*}

 
\subsection*{Sequences and bubble spaces}
Here we give the definitions of an exact sequence
and 
its sequence of traces. To emphasize that functions in a given finite dimensional space are defined on 
a domain of $\mathbb{R}^d$, we use the subscript $d$.

\begin{definition}[exact sequences] The sequences
\begin{alignat*}{3}
0\longrightarrow&\bR
\overset{i}{\longrightarrow}
H_3(K)
\overset{\grads}{\longrightarrow} 
E_3(K)
\overset{\curls}{\longrightarrow} 
V_3(K)
\overset{\divs}{\longrightarrow} 
W_3(K)
\overset{o}{\longrightarrow} 0
&&\quad\mbox{ for }K\subset \mathbb{R}^3,
\\
0\longrightarrow&\bR
\overset{i}{\longrightarrow}
H_2(K)
\overset{\grads}{\longrightarrow} 
E_2(K)
\overset{\curls}{\longrightarrow} 
W_2(K)
\overset{o}{\longrightarrow} 0
&&\quad\mbox{ for }K\subset \mathbb{R}^2,
\\
0\longrightarrow&\bR
\overset{i}{\longrightarrow}
H_1(K)
\overset{\grads}{\longrightarrow} 
W_1(K)
\overset{o}{\longrightarrow} 0
&&\quad\mbox{ for }K\subset \mathbb{R}^1,
\end{alignat*}
are said to be {\em exact} if

\

{\centering
\begin{tabular}{c c c}
$\mbox{ for } K\subset \mathbb{R}^3$
&$\mbox{ for } K\subset \mathbb{R}^2$
&$\mbox{ for } K\subset \mathbb{R}^1$
\\
\hline
\noalign{\smallskip}
\noalign{\smallskip}
$i(\bR)=\kerl_{\grads} H_3(K)$
&
$i(\bR)=\kerl_{\grads} H_2(K)$
&
$i(\bR)=\kerl_{\grads} H_1(K)$
\\
$\grads\! H_3(K)=\kerl_{\curls} E_3(K)$
&
$\grads \!H_2(K)=\kerl_{\curls} E_2(K)$
&
$\grads \!H_1(K)=\kerl_{o} W_1(K)$
\\
$\curls E_3(K)=\kerl_{\divs} E_3(K)$
&
$\curls E_2(K)=\kerl_{o} W_2(K)$
&
\\
$\divs V_3(K)=\kerl_{o} W_3(K)$
&
&
\end{tabular}
}
 \end{definition}


\begin{definition}[sequence of traces] Let $K\subset\mathbb{R}^3$ be a polyhedron. For any sequence
\begin{alignat*}{1}
\mathrm{S}(K):\quad 
0\longrightarrow\bR
\overset{i}{\longrightarrow}
H_3(K)
\overset{\grads}{\longrightarrow} 
E_3(K)
\overset{\curls}{\longrightarrow} 
V_3(K)
\overset{\divs}{\longrightarrow} 
W_3(K)
\overset{o}{\longrightarrow} 0,
\end{alignat*}
and, for any face $\fc\in \FC(K)$, we define its sequence of  traces on $\fc$ as  
\begin{alignat*}{2}
\tr^\fc\left(\mathrm{S(K)}\right)\mathrm{:} 
\quad 
0\longrightarrow\bR
\overset{i}{\longrightarrow}
\tr^\fc_H (H_3(K))
\overset{\grads}{\longrightarrow} 
\tr^\fc_E (E_3(K))
\overset{\curls}{\longrightarrow} 
\tr^\fc_V (V_3(K))
\overset{o}{\longrightarrow} 0,
\end{alignat*}
and, for any edge $\eg\in \EG(K)$, we define its sequence of  traces on $\eg$ by
\begin{alignat*}{2}
\tr^\eg\left(\mathrm{S(K)}\right)\mathrm{:} 
\quad 
0\longrightarrow\bR
\overset{i}{\longrightarrow}
\tr^\eg_H (H_3(K))
\overset{\grads}{\longrightarrow} 
\tr^\eg_E (E_3(K))
\overset{o}{\longrightarrow} 0,
\end{alignat*}
Sequence of traces on an edge for a two-dimensional sequence is defined in the same way.
\end{definition}

\begin{definition}[Bubble spaces] Let $K\subset\mathbb{R}^3$ be a polyhedron. For any sequence
\begin{alignat*}{1}
\mathrm{S}(K):\quad 
\bR
\overset{i}{\longrightarrow}
H_3(K)
\overset{\grads}{\longrightarrow} 
E_3(K)
\overset{\curls}{\longrightarrow} 
V_3(K)
\overset{\divs}{\longrightarrow} 
W_3(K)
\overset{o}{\longrightarrow} 0,
\end{alignat*}
we define the related $H^1$-, $H(\mathrm{curl})$-, $H(\mathrm{div})$-,
and $L^2$-{\em bubble} spaces as 
\begin{alignat*}{1}
\overset{\circ}{H_3}(K)&:=\{v\in H_3(K):\;\tr_Hv=0 \},
\\
\overset{\circ}{E_3}(K)&:=\{v\in E_3(K):\;\tr_Ev=0 \},
\\
\overset{\circ}{V_3}(K)&:=\{v\in V_3(K):\;\tr_Vv=0 \},
\\
\overset{\circ}{W_3}(K)&:=\{v\in W_3(K):\;\int_K v=0 \},
\end{alignat*}
respectively.

Similar, obvious definitions for bubble spaces hold for the two- and one-dimensional cases.
\end{definition}

From now on, we remove the first two and last terms 
in the definition of the exact sequence to simplify the notation. 
For example, we simply write
$H_1(\Omega)\overset{\grads}\longrightarrow W_1(\Omega)$, 
instead of writing
$0\longrightarrow\bR
\overset{i}{\longrightarrow} 
H_1(\Omega)\overset{\grads}\longrightarrow W_1(\Omega)\overset{o}{\longrightarrow} 0$.

Finally, let us emphasize that {\em all} the spaces in the sequences considered below
will have {\em finite} dimension.

\subsection*{Polynomial spaces}
We denote  the polynomial space of degree at most $p$ with argument $(x,y,z)\in \bR^3$ by
\[
 \pol_p(x,y,z) := \mathrm{span}\{x^iy^jz^k:\;i, j, k\ge 0,\, i+j+k\le p\},
\]
and we denote the {\it homogeneous} polynomial space of total degree $p$ by
\[
 \widetilde\pol_p(x,y,z) := \mathrm{span}\{x^iy^jz^k:\;i, j,k\ge 0,\, i+j+k= p\}.
\]
We denote the tensor-product polynomial space of degree at most $p$
by
\[
 \qol_p(x,y,z) := \pol_p(x)\otimes \pol_p(y)\otimes\pol_p(z) = 
 \mathrm{span}\{x^iy^jz^k:\;0\le i,j,k\le p\}
\]
We also denote the polynomial space of degree at most $p$ in the $(x,y)$ variable and 
of degree at most $p$ in the $z$ variable
by
\[
 \pol_{p|p}(x,y,z) := \pol_p(x,y)\otimes \pol_p(z).
\]
Similar definitions hold in the two- and one-dimensional cases.

Given an element $K\subset \bR^d$, we denote 
$\pol_p(K)$ to be the space of polynomials 
with degree at most $p$ defined on $K$. Similarly for 
$\widetilde\pol_p(K)$,
$\qol_p(K)$, $\pol_{p|p}(K)$.
We denote by $\bpol_p(K)$, respectively, $\widetilde\bpol_p(K)$, $\bqol_p(K)$, and 
$\bpol_{p|p}(K)$, the vector-valued functions whose components lie in 
$\pol_p(K)$, respectively, $\widetilde\pol_p(K)$,
 $\qol_p(K)$, and 
$\pol_{p|p}(K)$.

Whenever there is no possible confusion, we write $\pol_p$ instead of $\pol_p(K)$.

\subsection{Compatible exact sequences}
Here we introduce the concept of {\em compatible exact sequences}
in one-, two- and  three-space dimensions,
which is just a reformulation, in our notation, of the  {\em compatible} FES in differential form language introduced in \cite[Definition 5.12]{ChristiansenMuntheKaasOwren11}.

The main result of a compatible FES in \cite[Proposition 5.44]{ChristiansenMuntheKaasOwren11}, see also \cite[Proposition 2.8]{ChristiansenRapetti16}, provides the equivalence of a compatible FES with a FES admitting a commuting diagram. This powerful result reduces the search for a commuting diagram to that of a compatible FES (or compatible  exact sequence in our notation). 
In it, the {\em harmonic interpolator}, a generalization of the 
projection-based interpolation operator proposed  in \cite{DemkowiczBabuska03,DemkowiczBuffa05}, was used to obtain the commuting diagram; we reformulate these harmonic interpolators in one-, two- and  three-space dimensions using our notation in the Appendix A.


For a comprehensive theory of FESs, we refer to \cite[Section 5]{ChristiansenMuntheKaasOwren11} and \cite{ChristiansenRapetti16}.  See also \cite{ChristiansenGillette15} where the concept of a minimal compatible FES was introduced.


\begin{definition}[One-dimensional
compatible exact sequence]
\label{1d-cex}
Let $K$ be a segment.
Consider the finite dimensional exact sequence
\[
S_1(K):\quad H_1(K)\overset{\grads}{\longrightarrow} W_1(K).
\]
We say that the sequence $S_1(K)$ is a {\em compatible exact sequence} if
\begin{itemize}
\item [{\rm (i)}]
$ \dim \tr_H H_1(K) = 
\underset{\vt\in\VT(K)}{\sum} 1  =2$.
\end{itemize}
\end{definition}

\begin{definition}[Two-dimensional 
compatible exact sequence]
\label{2d-cex}
Let $K$ be a polygon.
Consider the finite dimensional exact sequence
\[
S_2(K):\quad H_2(K)\overset{\grads}{\longrightarrow}
 E_2(K)\overset{\curls}{\longrightarrow}
 W_2(K),
 \]
 and, for every edge $\eg\in \EG(K)$,  its sequence of traces
\[
\tr^\eg(S_2(K)):\quad H_1(\eg)\overset{\grads}{\longrightarrow}
 W_1(\eg),
\]
where $H_1(\eg):=\tr_H^\eg H_2(K)$ and
$W_1(\eg):=\tr_E^\eg E_2(K)$.
We say that the sequence $S_2(K)$ is a {\em compatible exact sequence} if
\begin{itemize}
\item [{\rm (i)}] For each edge $\eg\in\EG(K)$, 
the sequence $\tr^\eg(S_2(K))$
is a (one-dimensional) compatible exact sequence,
\item [{\rm (ii)}]  
$\begin{cases}\dim \tr_H(H_2(K))=\underset{\vt\in\mathcal{V}(K)}{\sum} 1
+\underset{\eg\in \EG(K)}{\sum}\dim \overset{\circ}{H_1}(\eg),\\
\dim \tr_E (E_2(K))=  
\underset{\eg\in \EG(K)}{\sum}\dim {W}_1(\eg),
\end{cases}$
\end{itemize}
\end{definition}

\begin{definition}[Three-dimensional
compatible exact sequence]
\label{3d-cex}
Let $K$ be a polyhedron.
Consider the  exact sequence
\[
S_3(K):\quad H_3(K)\overset{\grads}{\longrightarrow}
 E_3(K)\overset{\curls}{\longrightarrow}
 V_3(K)\overset{\divs}{\longrightarrow}
 W_3(K)
\]
and its sequences of traces for all faces $\fc\in\FC(K)$
\[
\tr^\fc(S_3(K)):\quad 
H_2(\fc)\overset{\grads}{\longrightarrow}
E_2(\fc)\overset{\curls}{\longrightarrow}
  W_2(\fc)
\]
and all edges $\eg\in\EG(K)$
 \[
\tr^\eg(S_3(K)):\quad H_1(\eg)\overset{\grads}{\longrightarrow}
 W_1(\eg),
\]
where $H_2(\fc)\times
E_2(\fc)\times
  W_2(\fc) := 
 \tr_H^\fc H_3(K)\times
 \tr_E^\fc E_3(K)\times
  \tr_V^\fc V_3(K),
  $ and 
  $H_1(\eg)\times
W_1(\eg) := 
 \tr_H^\eg H_3(K)\times
 \tr_E^\eg E_3(K).
  $ 

We say that the sequence $S_3(K)$ is a {\em compatible exact sequence} if
\begin{itemize}
\item [{\rm (i)}] For each face $\fc$, the sequence $\tr^\fc(S_3(K))$ 
is a (two-dimensional) compatible exact sequence.
\item [{\rm (ii)}] $\begin{cases} \dim \tr_H H_3(K) = \underset{\vt\in\mathcal{V}(K)}{\sum} 1
+\!\!\!\!\underset{\eg\in \EG(K)}{\sum}\!\!\!\!\dim \oo{H_1}(\eg) 
+\!\!\!\! \underset{\fc\in\FC(K)}{\sum}\!\!\!\!\dim \oo{H_2}(\fc),
\\
\dim \tr_E E_3(K)= \underset{\eg\in \EG(K)}{\sum} W_1(\eg)
+\underset{\fc\in\mathcal{F}(K)}{\sum}\dim \overset{\circ}{E}_2(\fc)
, \\
\dim \tr_V V_3(K)=\underset{\fc\in\FC(K)}{\sum}
               \dim W_2(\fc).
                           \end{cases}$
\end{itemize}
\end{definition}

\

The following result on the equivalence of compatible exact sequence and compatible FES is trivial to 
verify by definition; we omit its proof.
\begin{proposition}
Let $K\in \bR^d$ ($d=1,2,3$) be a $d$-dimensional polytope. Then, 
with a change of scalar/vector fields to differential forms,
a compatible exact sequence on $K$ defined above
is the corresponding compatible FES on $\triangleleft{(K)}$, 
the set that contains the element $K$ and its vertices, edges (if $d\ge 2$), and faces (if $d= 3$).
\end{proposition}

\

The following result is a direct consequence of  \cite[Proposition 2.8]{ChristiansenRapetti16}.
\begin{proposition}
Let $K\in \bR^d$ ($d=1,2,3$) be a $d$-dimensional polytope. The following statements are equivalent:
\begin{itemize}
 \item $S(K)$ is a compatible exact sequence on $K$.
 \item $S(K)$ admits a commuting diagram.
\end{itemize}
\end{proposition}
Note that, in the proof of \cite[Proposition 2.8]{ChristiansenRapetti16}, the so-called 
{\em harmonic interpolator} is used to obtain the commuting diagram. We give a reformulation of this concept in our notation in the Appendix A.

\subsection{The construction of compatible exact sequences}
Here, we give our main results on the systematic construction of compatible exact sequences. Again, for each space dimension, we have a specific construction. The corresponding proofs are given in Section \ref{sec:proof}. The resulting compatible exact sequences are all minimal \gfu{\jbc{sequences containing} a prescribed exact sequence}
in the sense of 
\cite[Corollary 3.2]{ChristiansenGillette15}.

\subsubsection*{The one-dimensional case}
In one space dimension, the construction of compatible exact sequences is fairly simple.
\begin{theorem}\label{thm:1d-general}
 Let $K$ be a 
 segment.
 Let any given  exact sequence
\[
S^g_1(K):\quad H^g_1(K)\overset{\grads}{\longrightarrow} W^g_1(K)
\]
be such that
\[
\pol_0(K)\subset 
W_1^g(K).
\]
 Then, it 
 is {compatible}.
\end{theorem}

\subsubsection*{The two-dimensional case}
In two-space dimensions, the systematic construction of compatible exact sequence is more involved than in the one-dimensional case, not only because of the geometry but also because we seek a 
compatible exact sequence with certain given traces on each of the edges of the element. 
Those traces are compatible exact sequences (we found while dealing with the one-dimensional case) 
which we gather in the set $S_1(\dK)$. 
The sequence with which we begin the construction must then be what we call 
{\em $S_1(\dK)$-admissible}. We define this term next.

\begin{definition}[$S_1(\dK)$-admissible exact sequence] 
\label{s1-ad}
Let 
\[
S_1(\dK):=\{S_1(\eg):\quad
H_1(\eg)\overset{\grads}{\longrightarrow} W_1(\eg)\quad\forall e\in \EG(K)\},
\]
be any set of 
one-dimensional
exact sequences. 
Then,
we say that a given 
two-dimensional  exact sequence
\[
 H_2^g(K)\overset{\grads}{\longrightarrow} E_2^g(K)
 \overset{\curls}{\longrightarrow} W_2^g(K)
\]
is $S_1(\dK)$-admissible if
\begin{itemize}
\item [{\rm (i)}]
$\tr_H^\eg\left( H_2^g(K)\right)\times
 \tr_E^\eg\left( E_2^g(K)\right)
 \subset H_1(\eg)\times W_1(\eg)\quad\forall \mbox{ e }\in \EG(K),$
\item[{\rm (ii)}] $\pol_0(K)\subset W_2^g(K)$.
\end{itemize}
\end{definition}

The next theorem is the main result of the two-dimensional case. It shows how to construct a compatible exact sequence by suitably enriching a given $S_1(\dK)$-admissible exact sequence. The compatible exact sequence we seek is such that its traces on the edges {\em coincide} with the exact sequences of the set $S_1(\dK)$.

\begin{theorem}\label{thm:2d-general}
Let $K$ be a polygon, let 
\[
S_1(\dK):=\{S_1(\eg):\quad
H_1(\eg)\overset{\grads}{\longrightarrow} W_1(\eg)\quad\forall e\in \EG(K)\},
\]
be a set of {\em compatible exact sequences}, 
 and let
 \[
 H_2^g(K)\overset{\grads}{\longrightarrow} E_2^g(K)
 \overset{\curls}{\longrightarrow} W_2^g(K)
\]
be a {\em given} $S_1(\dK)$-admissible  exact sequence. 
Let the space $\delta H_2^g(K)\subset H^1(K)$
 satisfy the following properties:
\begin{itemize}
 \item [(i)] $\tr_H^\eg\delta H_2^g(K)\subset H_1(\eg) $ for all edges $\eg\in \EG(K)$.
 \item [(ii)] $\delta H_2^g(K)\cap  H_2^g(K) =\{0\}$.
 \item [(iii)] $
 \{v\in H_2^g(K)\oplus\delta H_2^g(K):\;\tr_H v=0\}
 =\oo{H_2^g}(K).$
 \item [(iv)] $\dim \delta H_2^g(K) = \underset{\vt\in \VT(K)}{\sum} 1 +
 \underset{\eg\in \EG(K)}{\sum} \dim \oo{H_1}(e) + \dim 
 \oo{H_2^g}(K)-\dim {H_2^g}(K)$.
\end{itemize}
Then, the following sequence 
\begin{equation*}
\begin{tabular}{c  c c c c}
$S_2(K):\quad{H_2^g(K)\oplus \delta H_2^g(K)}$&
$\overset{\grads}{\longrightarrow}$ &
$E_2^g(K)\oplus\grads\delta H_2^g(K)$
&
$\overset{\curls}{\longrightarrow}$ &
$W_2^g(K)$\\
\end{tabular}
\end{equation*}
is a compatible exact sequence.
Moreover, 
it is also a minimal compatible exact sequence
containing the exact sequence
 \[
 H_2^g(K)\overset{\grads}{\longrightarrow} E_2^g(K)
 \overset{\curls}{\longrightarrow} W_2^g(K).
\]
\end{theorem}

\

Let us relate this result with the theory of $M$-decompositions 
developed in \cite{CockburnFuSayas16}. Such theory has to do with the right-most part of the commuting diagram.
Since the operator $\nabla\times$ was replaced by the divergence operator $\nabla\cdot$ and since
\[
\nabla\times (v_1,v_2)=-\partial_{x_2} v_1 +\partial_{x_1} v_2=\nabla\cdot (v_2, -v_1)=\nabla\cdot (v_1,v_2)^{\mathrm{rot}},
\]
we have that $\nabla\times E_2(K)=\nabla\cdot E_2^{\mathrm{rot}}(K)$, with the obvious notation. 
In \cite[Proposition 5.1]{CockburnFuSayas16} is was shown that
$\left(E_2^g(K)\oplus \grads \delta H_2^g(K)\right)^{\mathrm{rot}}\times W_2^g(K)$ is the {\it smallest}
space  containing $\left(E_2^g(K)\right)^{\mathrm{rot}}\times W_2^g(K)$ which admits an $M(\dK)$-decomposition with
the trace space 
\[
 M(\dK):=\{\mu\in L^2(\dK):\;\;
 \mu|_{\eg}\in W_1(\eg)\quad \forall \eg\in\EG(K)\}.
\]
This result implies that the  right-most part of the diagram commutes.

\subsubsection*{The three-dimensional case}
Now, we present our most involved result in three-space dimensions.
As we did in the two-dimensional case, we have to provide a set of compatible exact sequences on the faces of the element $K$, $S_2(\partial K)$;
they were previously obtained when dealing with the two-dimensional case.
Then, we have to  start from certain exact sequences we call 
{\em $S_2(\dK)$-admissible}. 

\begin{definition}[$S_2(\dK)$-admissible exact sequence]
\label{s2-ad}
Let 
\[
S_2(\dK):=\{S_2(\fc):\quad H_2(\fc)\overset{\nabla}{\longrightarrow}
E_2(\fc)\overset{\nabla\times}{\longrightarrow}W_2(\fc), 
\; \fc\in\mathcal{F}(K)\}
\]
be any set of two-dimensional  exact sequences.
We say that a given three-dimensional exact sequence
\[
 H_3^g(K)\overset{\grads}{\longrightarrow} E_3^g(K)
 \overset{\curls}{\longrightarrow} V_3^g(K)
\overset{\divs}{\longrightarrow} W_3^g(K)
\]
is $S_2(\dK)$-admissible if
\begin{itemize}
\item[{\rm(i)}] $\tr_H^\fc H_3(K)\times \tr_E^\fc E_3^g(K)
\times \tr_V^\fc V_3^g(K) \subset H_2(\fc)\times E_2(\fc)\times W_2(\fc)\;\forall\fc\in \mathcal{F}(K)$,
\item[{\rm (ii)}] $\pol_0(K)\subset W_3^g(K)$.
\end{itemize}
\end{definition}

The following theorem is our main result. It shows how to enrich an $S_2(\dK)$-admissible exact sequence to get a commuting exact sequence.
Note also that the traces on the faces  of the sequence we seek must {\em coincide} with the sequence of traces in the set $S_2(\dK)$. 
\begin{theorem}
\label{thm:3d-general}
Let $K$ be a polyhedron, let
\[
S_2(\dK):=\{S_2(\fc):\quad H_2(\fc)\overset{\nabla}{\longrightarrow}
E_2(\fc)\overset{\nabla\times}{\longrightarrow}W_2(\fc), 
\; \fc\in\mathcal{F}(K)\}
\]
be a set of {\em compatible exact sequences} satisfying the following {\em compatibility} condition
\[
 \tr_H^\eg H_2(\fc_1)\times
 \tr_E^\eg E_2(\fc_1)
=  \tr_H^\eg H_2(\fc_2)
\times
 \tr_E^\eg E_2(\fc_2)=:H_1(\eg)\times W_1(\eg),
\]
for any faces $\fc_1, \fc_2\in \FC(K)$ sharing an 
edge $\eg: =\EG(\fc_1)\cap \EG(\fc_2)$. 
Let
\[
 H_3^g(K)\overset{\grads}{\longrightarrow} E_3^g(K)
 \overset{\curls}{\longrightarrow} V_3^g(K)
\overset{\divs}{\longrightarrow} W_3^g(K)
\]
be a {\em given} $S_2(\dK)$-admissible  
exact sequence.

Let the spaces $\delta H_3^g(K)\times\delta E_3^g(K) \subset H^1(K)\times 
 H(\mathrm{curl},K)
$
 satisfy the following properties:
 
 \
 
\underline{Properties of $\delta H_3^g(K)$} 

\

\begin{itemize}
 \item [(i)] $\tr_H^\fc\delta H_3^g(K)\subset H_2(\fc) $ for all faces  $\fc\in \FC(K)$.
 \item [(ii)] $\delta H_3^g(K)\cap  H_3^g(K) =\{0\}$.
 \item [(iii)] $
 \{v\in H_3^g(K)\oplus\delta H_3^g(K):\;\tr_H v=0\}
 =\oo{H_3^g}(K).$
 \item [(iv)] $\dim \delta H_3^g(K) =\underset{\vt\in \VT(K)}{\sum}1
+\underset{\eg\in \EG(K)}{\sum}\dim \oo{H_1}(\eg) 
+ \underset{\fc\in\FC(K)}{\sum}\dim \oo{H_2}(\fc)$\\
$\phantom{ooooooooooooooo}+\dim \oo{H_3^g}(K)-\dim {H_3^g}(K)$.
\end{itemize}

\
 
\underline{Properties of $\delta E_3^g(K)$} 

\

\begin{itemize}
 \item [(i)] $\tr_E^\fc\delta E_3^g(K)\subset E_2(\fc) $ for all faces  $\fc\in \FC(K)$.
 \item [(ii)] $\curls\delta E_3^g(K)\cap  V_3^g(K) =\{0\}$.
 \item [(iii)] $
 \{v\in V_3^g(K)\oplus\curls\delta E_3^g(K):\;\tr_V v=0,\;
 \divs v=0\}
 =\\
\{v\in \oo{V_3^g}(K):\;
 \divs v=0\}.$
 \item [(iv)] 
               $\dim \delta E_3^g(K) = 
               \dim \curls\delta E_3^g(K) =\\
               \underset{\fc\in\FC(K)}{\sum}
               \dim W_2(\fc) +\dim\oo{W_3^g}(K) +\dim \{v\in \oo{V_3^g}(K):
               \;\divs v=0\}-\dim {V_3^g}(K)$.
\end{itemize}

Then, the sequence
\begin{equation*}
\resizebox{01.\columnwidth}{!}{
\begin{tabular}{c c c c c c c}
$S_3(K):\quad{H_3^g(K)}$&
$\overset{\grads}{\longrightarrow}$ &
$E_3^g(K)$
&
$\overset{\curls}{\longrightarrow}$ &
$V_3^g(K)$
&
$\overset{\divs}{\longrightarrow}$ &
$W_3^g(K)$
\\
$\phantom{S_3(K):\quad}\oplus\delta H_3^g(K)$&
&
$\oplus\grads\delta H_3^g(K)
\oplus\delta E_3^g(K)$
&
&
$\oplus\curls \delta E_3^g(K)$
\\
\end{tabular}}
\end{equation*}
is a commuting exact sequence.
Moreover, it is a minimal commuting exact sequence
containing the exact sequence
\[
 H_3^g(K)\overset{\grads}{\longrightarrow} E_3^g(K)
 \overset{\curls}{\longrightarrow} V_3^g(K)
\overset{\divs}{\longrightarrow} W_3^g(K).
\]
\end{theorem}

\

Note that, unlike the two-dimensional case, we are requiring the sequences of the set $S_2(\dK)$ to satisfy a compatibility condition on each of the edges of the polyhedral. Such compatibility condition was automatically satisfied, and hence was not required, in the two-dimensional case.

Let us relate this result with the theory of $M$-decompositions 
introduced in \cite{CockburnFuSayas16}. As for the two-dimensional case, the part of this theory concerned with mixed methods is associated 
with the right-most side of the diagram. Indeed, in \cite[Proposition 5.1]{CockburnFuSayas16} it is shown that the space
$\left(V_3^g(K)\oplus \curls \delta E_3^g(K)\right)\times W_3^g(K)$ is the {\it smallest}
one containing $V_3^g(K)\times W_3^g(K)$ and admitting an $M(\dK)$-decomposition with
the trace space 
\[
 M(\dK):=\{\mu\in L^2(\dK):\;\;
 \mu|_{\fc}\in W_2(\fc)\quad \forall \fc\in\FC(K)\}.
\]

\section{Applications}
\label{sec:applications}
Now, we apply our main results on the systematic construction of compatible  exact sequences in 
Section \ref{sec:mainresults} to explicit  construct them on various element shapes in one-, two- and  three-space dimensions. 
To emphasis on the existence of a commuting diagram for compatible  exact sequence, we denote such sequence as {\em commuting exact sequence}.

\subsection{The one-dimensional case}
In one-space dimension, we consider the element 
$K\subset\bR^1$ to be the reference 
interval $\{x:\;0<x<1\}$.
The most useful commuting exact sequence on $K$ is given below.
\begin{theorem}\label{thm:1d}
Let $K\subset \bR^1$ be the reference interval 
with coordinate $x$. Then,
the following sequence on $K$ is a commuting exact sequence
for $k\ge 0$.
\begin{alignat*}{2}
\pol_{k+1}(x)
\overset{\grads}{\longrightarrow} 
\pol_{k}(x)
\end{alignat*}
\end{theorem}

\subsection{The two-dimensional case}
In two-space dimensions, we proceed as follows. 
We first consider the element $K$ to be  
either the reference triangle $\{(x,y):\;x>0,y>0, x+y<1 \}$, or 
the reference square $\{(x,y):\;0<x<1, 0<y<1\}$.
We present two commuting exact sequences 
on the reference triangle 
 and four on the reference square. 
All the sequences, except the second one on the reference square, are known and its spaces are all spaces of polynomials. 

Then, we consider the case in which $K$ is a general polygon, and present two {\it new} commuting exact sequences which contain 
non-polynomial functions; this is based on the results in  \cite{CockburnFuM2D} on the construction of $M$-decompositions in two dimensions.
All these commuting exact sequences are the {\it smallest} ones, as stated in Theorem \ref{thm:2d-general},
that contain certain given  exact sequence and have 
a certain prescribed sequence of traces on each edge.

We end by obtaining commuting exact sequences on two-dimensional polygonal meshes. 

\subsection*{Triangle}

\

\begin{theorem}\label{thm:2d-t}
Let $K$ be the reference triangle with coordinates $(x,y)$. 
Then,
the following two sequences on $K$ are  
commuting exact sequences
for $k\ge 0$,

\resizebox{0.95\columnwidth}{!}{
\begin{tabular}{l c c c c c}
$
\satri(K) :
$ &
$ \pol_{k+2}(x,y)$&
$\overset{\grads}{\longrightarrow}$ &
$\bld\pol_{k+1}(x,y)$&
$\overset{\curls}{\longrightarrow}$ &
$\pol_{k}(x,y)$,\\
$
\sbtri(K) :$ &$\pol_{k+1}(x,y)$&
$\overset{\grads}{\longrightarrow}$ &
$\bld\pol_{k}(x,y)\oplus\bld x
\times\widetilde\pol_{k}(x,y)$
&
$\overset{\curls}{\longrightarrow}$ &
$\pol_{k}(x,y)$.
\end{tabular}}\\
Here $\bld x\times p = (y\,p,-x\,p)^t$ for a scalar function $p$.

Moreover, 
the sequence of traces for $\satri$ on an edge $\eg\in\EG(K)$
is 
\begin{subequations}
\begin{align*}
 \tr^\eg\left(\satri(K)\right): 
 \pol_{k+2}(\eg)\longrightarrow \pol_{k+1}(\eg), 
\end{align*}
and that for $\sbtri$ is 
\begin{align*}
 \tr^\eg\left(\sbtri(K)\right): 
 \pol_{k+1}(\eg)\longrightarrow \pol_{k}(\eg). 
\end{align*}
\end{subequations}
\end{theorem}

\

{
These two sequences are well-known. Indeed,
the first sequence $\satri$ is mainly due-to Brezzi, 
Douglas, and Marini 
\cite{BrezziDouglasMarini85}
since its $H(\mathrm{curl})$ space is a ninety-degree rotation of 
the $H(\mathrm{div})$-space, usually called the BDM space, obtained in \cite{BrezziDouglasMarini85},
of degree $k+1$. Its $H^1$ and $L^2$ spaces are the Lagrange polynomial spaces of 
degree $k+2$ and discontinuous polynomial
space of degree $k$.
The second sequence $\sbtri$ is mainly due-to 
Raviart and Thomas 
\cite{RaviartThomas77}
since its $H(\mathrm{curl})$ space is a ninety-degree rotation of 
the $H(\mathrm{div})$ space, usually called the RT space,  obtained in \cite{RaviartThomas77}, 
of degree $k$. Its $H^1$ and $L^2$ spaces are the Lagrange polynomial spaces of 
degree $k+1$ and discontinuous polynomial
space of degree $k$.
}

\subsection*{Square}

\

\begin{theorem}\label{thm:2d-s}
Let $K$ be the reference square with coordinates $(x,y)$. Then,
the following four sequences are commuting exact sequences 
for $k\ge 0$:

\resizebox{0.95\columnwidth}{!}
{
\begin{tabular}{l c c c c c}
$\sasqr(K) :$ &
$ \pol_{k+2}(x,y)$&
$\overset{\grads}{\longrightarrow}$ &
$\bld\pol_{k+1}(x,y) $&
$\overset{\curls}{\longrightarrow}$ &
$\pol_{k}(x,y)$,\\
& $\oplus\delta H^{2,I}_{k+2}$ & 
& $\oplus\grads\delta H^{2,I}_{k+2}$\\
$\sbsqr(K):$ &$\pol_{k+1}(x,y)$&
$\overset{\grads}{\longrightarrow}$ &
$\bld\pol_{k}(x,y)\oplus\bld x
\times\widetilde\pol_{k}(x,y)$
&
$\overset{\curls}{\longrightarrow}$ &
$\pol_{k}(x,y)$,\\
& $\oplus\delta H^{2,I}_{k+1}$ & 
& $\oplus\grads\delta H^{2,I}_{k+1}$
\vspace{.3cm}\\
\hline
\\
$\scsqr(K):$ &$\qol_{k}(x,y)\oplus \{x^{k+1},y^{k+1}\}$&
$\overset{\grads}{\longrightarrow}$ &
$\bld\qol_{k}(x,y)\oplus\bld x
\times\{x^ky^k\}$
&
$\overset{\curls}{\longrightarrow}$ &
$\qol_{k}(x,y)$,\\
& $\oplus\delta H^{2,I}_{k+1}$ & 
& $\oplus\grads\delta H^{2,I}_{k+1}$\\
$\sdsqr(K):$ &$\qol_{k+1}(x,y)$&
$\overset{\grads}{\longrightarrow}$ &
$\bld\qol_{k}(x,y)\oplus
\left(
\begin{tabular}{c}
 $y^{k+1}\pol_k(x)$\\
 $x^{k+1}\pol_k(y)$
\end{tabular}
\right)$
&
$\overset{\curls}{\longrightarrow}$ &
$\qol_{k}(x,y)$.
\end{tabular}}\\
Here the additional space $\delta H^{2,I}_{k}$, for $k\ge 1$,
takes the following form: 
\begin{alignat*}{2}
 \delta H^{2,I}_{k} :=&\;\mathrm{span} \{x\,y^{k},y\,x^k\}.
\end{alignat*}
Moreover, 
the sequence of traces for $\sasqr$ on an 
edge $\eg\in\EG(K)$ is 
\begin{subequations}
\begin{align*}
 \tr^\eg\left(\sasqr(K)\right): 
 \pol_{k+2}(\eg)\longrightarrow \pol_{k+1}(\eg), 
\end{align*}
and that for $\sisqr$  with $i\in\{2,3,4\}$ is 
\begin{align*}
 \tr^\eg\left(\sisqr(K)\right): 
 \pol_{k+1}(\eg)\longrightarrow \pol_{k}(\eg). 
\end{align*}
\end{subequations}
\end{theorem}

\

Note that for $k=0$, the last three sequences are 
exactly the same.
{
Here the second sequence is new and the other three are 
well-known. 
The first sequence $\sasqr$ is mainly due-to 
Brezzi, Douglas, and Marini \cite{BrezziDouglasMarini85}
since its $H(\mathrm{curl})$ space is a ninety-degree rotation of 
the $H(\mathrm{div})$ space, usually called the BDM space, obtained in  \cite{BrezziDouglasMarini85}, 
of degree $k+1$ on the square. 
Its $H^1$ and $L^2$ spaces are the serendipity polynomial spaces of 
degree $k+2$ and discontinuous polynomial
space of degree $k$.
The second sequence $\sbsqr$ is a new one resulting a new family of 
$H(\mathrm{curl})$ spaces. 
The third sequence is the TNT sequence \cite{CockburnQiu14} on the square.
And the last one is mainly due to Raviart and Thomas \cite{RaviartThomas77}
since its $H(\mathrm{curl})$ space is a ninety-degree rotation of 
the $H(\mathrm{div})$ space, usually called the RT space, obtained in \cite{RaviartThomas77},
of degree $k$ on the square. 
Its $H^1$ and $L^2$ spaces are the tensor-product Lagrange polynomial space of 
degree $k+1$ and discontinuous tensor product polynomial
space of degree $k$.
}

\subsection*{Polygon}

\

The explicit construction 
of (high-order) commuting exact sequences 
on a general polygon $K$ is not known. 
Here we fill this gap by presenting two families of 
commuting exact sequences 
by applying Theorem \ref{thm:2d-general}. To do so, we take advantage of
the recent results on constructing
$M$-decompositions 
in \cite{CockburnFuM2D} to deal with the right-most is of the diagram.

To state our result, we need to introduce some notation. 
Let $\{\vt_i\}_{i=1}^{ne}$ be the set of vertices of the 
polygonal element $K$ which we take to be counter-clockwise 
ordered. 
Let $\{\eg_i\}_{i=1}^{ne}$ be the set of edges 
of $K$ where the edge $\eg_i$ connects the vertices 
$\vt_i$ and $\vt_{i+1}$.
Here the subindexes are integers module  $ne$.
We also define, for $1\le i\le ne$,  $\lambda_i$ to be the 
linear function that 
vanishes on edge $\eg_i$ and reaches maximum value $1$ 
in the closure of the  element $K$.
To each vertex $\vt_i$, $i=1,\dots,ne$, we associate a function
$\xi_i$ satisfying the following conditions:
\begin{itemize}
\item[(L.1)] $\xi_i \in H^1(K)$,
\item[(L.2)] $\xi_i|_{\eg_j}\in \pol_1(\eg_j),\;j=1,\dots,ne,$
\item[(L.3)] $\xi_i(\vt_j)=\delta_{i,j},\;j=1,\dots,ne,$
\end{itemize}
where $\delta_{i,j}$ is the Kronecker delta.
Note that conditions (L.2) and (L.3) together ensure 
that the trace of $\xi_i$ on the edges is only no-zero at 
$\eg_{i}$ and $\eg_{i+1}$, where they are linear.
These functions, with examples given in \cite{CockburnFuM2D},
 are not polynomials if the polygon $K$ is 
not a triangle or a parallelogram.

Now, we are ready to state the result.
\begin{theorem}
\label{thm:2d-p}
Let $K$ be a polygon of $ne$ edges
with coordinates $(x,y)$ that does not
have hanging nodes.
Then,
the following two sequences on $K$ are  
commuting exact sequences
for $k\ge 0$,

{
\begin{tabular}{l c c c c c}
$
\mathrm{S_{1,k}^{poly}}(K) :
$ &
$ \pol_{k+2}(x,y)$&
$\overset{\grads}{\longrightarrow}$ &
$\bld\pol_{k+1}(x,y)$&
$\overset{\curls}{\longrightarrow}$ &
$\pol_{k}(x,y)$,\\
&$ \oplus\delta H^{2,II}_{k+2}$&
 &
$\oplus\grads \delta H^{2,II}_{k+2}$&
 &\\
$
\mathrm{S_{2,k}^{poly}}(K) :$ &$\pol_{k+1}(x,y)$&
$\overset{\grads}{\longrightarrow}$ &
$\bld\pol_{k}(x,y)\oplus\bld x
\times\widetilde\pol_{k}(x,y)$
&
$\overset{\curls}{\longrightarrow}$ &
$\pol_{k}(x,y)$\\
&$\oplus \delta H^{2,II}_{k+1}$&
 &
$\oplus\grads \delta H^{2,II}_{k+1}$&
 &
\end{tabular}}\\
Here the additional space $\delta H_k^{2,II}$, for $k\ge 1$, takes 
the following form:
\begin{align*}
 \delta H_k^{2,II} = \oplus_{i=3}^{ne} \Psi_{i,k},
\end{align*}
where 
 \begin{alignat*}{2}
  \Psi_{i,k} = &\;\left\{ 
  \begin{tabular}{l c}
  $ \mathrm{span}\{\xi_{i+1}\lambda_{i+1}^a :\;\; 
 \max\{k+3-i,0\}\le a \le k-1\}$ &$\text{ if $3\le i \le ne-1$}$,\\
 $\mathrm{span}\{\xi_{i+1}\lambda_{i+1}^a :\;\; 
 \max\{k+4-i,1\}\le a \le k-1\}$ &$\text{ if $i = ne$}$,
 \end{tabular}
 \right.
 \end{alignat*}
and the functions $\{\xi_i\}_{i=1}^{ne}$ are assumed to 
satisfy conditions (L).

Moreover, 
the sequence of traces for $\mathrm{S_{1,k}^{poly}}$ on an edge $\eg\in\EG(K)$
is 
\begin{subequations}
\begin{align*}
 \tr^\eg\left(\satri(K)\right): 
 \pol_{k+2}(\eg)\longrightarrow \pol_{k+1}(\eg), 
\end{align*}
and that for $\mathrm{S_{2,k}^{poly}}$ is 
\begin{align*}
 \tr^\eg\left(\sbtri(K)\right): 
 \pol_{k+1}(\eg)\longrightarrow \pol_{k}(\eg). 
\end{align*}
\end{subequations}
\end{theorem}

\

{
Note that in the notation of \cite{CockburnFuM2D},
we have $\delta V_{\!\mathrm{fillM}}:=
\curls \delta H^{2,II}_{k+1}$
 is 
 the filling space to guarantee $M$-decompositions for 
 the pair 
 $V\times W:= \bpol_k\oplus \delta V_{\!\mathrm{fillM}} \times \pol_k$ for the trace space 
 \[
  M(\dK):=\{
  \mu\in L^2(\dK):
  \mu|_{\eg}\in \pol_k(\eg)\quad \forall \eg\in \EG(K)
  \}.
 \]
}

{When $K$ is a convex quadrilateral, construction of 
$H(\mathrm{div})$-conforming spaces of the form $ \bpol_k(K)\oplus \curls \delta H(K)$ 
was also presented in \cite{ArbogastCorrea15}. 
While the filling space $\curls \delta H(K)$ share the same dimension ($1$ if $k = 0$ and $2$ if $k\ge 1$) as those in \cite{CockburnFuM2D} on a convex quadrilateral, 
they were 
not constructed directly on the quadrilateral element $K$ as was done in \cite{CockburnFuM2D}
but by Piola mapping two (one if $k=0$) divergence-free polynomial functions
from the reference square to $K$, resulting rational functions. 
Moreover, $H(\mathrm{div})$-conforming finite element shape functions were also created 
\cite[Section 5]{ArbogastCorrea15}.
}

{
Let us now briefly comment on the special cases when
$K$ is a triangle or a parallelogram. 
 In these cases, the element $K$ can be considered as 
 a physical element obtained from 
 an affine mapping of the 
reference triangle or the reference square.
We can easily obtain {\it mapped} 
commuting exact sequences on the physical element $K$ 
from those 
in Theorem \ref{thm:2d-t} on the reference triangle, and 
those 
in Theorem \ref{thm:2d-s} on the reference square via 
proper linear mapping functions.
To simplify the notation, we still denote the mapped sequences on a physical
triangle $K$ as 
$\satri(K)$ and 
$\sbtri(K)$, and those on the physical 
parallelogram as $\sisqr(K)$ for $i\in\{1,2,3,4\}$.
These mapped sequences have an advantage over those in 
Theorem \ref{thm:2d-p}, namely, that numerical integration 
only needs to be done on the reference elements.
}

 \subsection*{Whole mesh} 
 Using the sequences in Theorem \ref{thm:2d-t}, and Theorem \ref{thm:2d-s}, 
 we can readily obtain {\it mapped} 
 commuting exact sequences on a hybrid mesh, $\Omega_h:=\{K\}$, of 
a polygonal domain $\Omega$, where each physical element $K$ is 
an affine mapping of the 
reference triangle or the reference square.
We obtain two families of {\it mapped } 
sequences with the following form,  
\[
\mathrm{S_{2d}}(\Omega_h):\quad  H_2(\Omega_h)\overset{\grads}{\longrightarrow} E_2(\Omega_h)
 \overset{\curls}{\longrightarrow} W_2(\Omega_h),
\]
where $H_2(\Omega_h)\times E_2(\Omega_h)
 \times W_2(\Omega_h)\subset 
 H^1(\Omega)\times H(\mathrm{curl},\Omega)\times L^2(\Omega)$.
The restriction on an element $K$ 
of 
the first family of sequences 
is $\satri(K)$ if $K$ is a triangle, and 
is $\sasqr(K)$ if $K$ is a parallelogram;
and 
the restriction on an element $K$ 
of 
the second family of sequences 
is $\sbtri(K)$ if $K$ is a triangle, and 
is any of  $\sisqr(K)$  for $i\in \{2,3,4\}$ 
if $K$ is a parallelogram.
The reason we are able to do so is due to 
the trace compatibility of the sequences $\satri(K)$
and $\sasqr(K)$, and 
the trace compatibility of the sequences $\sbtri(K)$
and $\sisqr$  for $i\in\{2,3,4\}$ in Theorem \ref{thm:2d-t} and 
Theorem \ref{thm:2d-s}. 

On the other hand, using the sequences in Theorem \ref{thm:2d-p}, 
 we can readily obtain
two families of {\it non-mapped} 
commuting exact sequences on a more general polygonal mesh, $\Omega_h:=\{K\}$, of 
a polygonal domain $\Omega$, where each physical element $K$ is 
a polygon, and 
the restriction on an element $K$ 
of the first family of sequences 
is $\mathrm{S_{1,k}^{poly}}(K)$, and that for 
the second 
is $\mathrm{S_{2,k}^{poly}}(K)$. 

\subsection{The three-dimensional case}
In three-space dimensions, we first consider the element $K\subset\bR^3$ to be  
either of the following four reference polyhedra:

\

\begin{tabular}{l l}
 ${
\begin{tikzpicture}
\draw  (0,0) -- (0,0.3);
\draw  (0,0) -- (0.2996,-0.0157);
\draw  (0,0) -- (-.2194,-0.2046);
\draw (0,0.3) -- (-.2194,-0.2046)--(0.2996,-0.0157)--(0,0.3);
\end{tikzpicture}}$
tetrahedron:
&
 $\{(x,y,z): 0<x, \;0<y,\; 0<z,\; x+y+z<1\}$ 
,\\
 ${
\begin{tikzpicture}
\draw  (0,0) -- (0,0.3);
\draw  (0,0) -- (0.2996,-0.0157);
\draw  (0,0) -- (-.2194,-0.2046);
\draw (-.2194,-0.2046) -- (.0802,-0.2046) -- (0.2996,-0.0157);
\draw  (0,0.3) -- (0.2996,0.2843);
\draw (0,0.3) -- (-.2194,0.0954);
\draw (-.2194,0.0954) -- (.0802,0.0954) -- (0.2996,0.2843);
\draw (-.2194,0.0954) -- (-.2194,-0.2046);
\draw (.0802,0.0954) -- (.0802,-0.2046);
\draw (0.2996,-0.0157) -- (0.2996,0.2843);
\end{tikzpicture}}$ 
hexahedron: &
 $\{(x,y,z): 0<x<1,\; 0<y<1, \;0<z<1\}$,\\
 ${
\begin{tikzpicture}
\draw  (0,0) -- (0,0.3);
\draw  (0,0) -- (0.2996,-0.0157);
\draw  (0,0) -- (-.2194,-0.2046);
\draw (-.2194,-0.2046) -- (0.2996,-0.0157);
\draw  (0,0.3) -- (0.2996,0.2843);
\draw (0,0.3) -- (-.2194,0.0954);
\draw (-.2194,0.0954)-- (0.2996,0.2843);
\draw (-.2194,0.0954) -- (-.2194,-0.2046);
\draw (0.2996,-0.0157) -- (0.2996,0.2843);
\end{tikzpicture}}$ 
prism: &
 $\{(x,y,z): 0<x,\; 0<y,\; 0<z<1,\; x+y<1\}$,\\
 ${
\begin{tikzpicture}
\draw  (0,0) -- (0,0.3);
\draw  (0,0) -- (0.2996,-0.0157);
\draw  (0,0) -- (-.2194,-0.2046);
\draw (-.2194,-0.2046) -- (.0802,-0.2046) -- (0.2996,-0.0157);
\draw (0,0.3) -- (-.2194,-0.2046);
\draw (0,0.3) -- (.0802,-0.2046);
\draw (0,0.3) -- (0.2996,-0.0157);
\end{tikzpicture}}$ 
pyramid: &
 $\{(x,y,z): 0<x,\; 0<y,\; 0<z,\; x+z<1,\;y+z<1\}$. 
\end{tabular}

\

We present two commuting exact sequences on the reference tetrahedron 
 and four on the reference hexahedron, prism, and pyramid. 
 All these commuting exact sequences are the {\it smallest} ones, as described in
 Theorem \ref{thm:3d-general}, which
contain certain given  exact sequence and have 
a certain prescribed sequence of traces on each face.
We then obtain commuting exact sequences
 on polyhedral meshes made of mapped tetrahedra
 hexahedra, prisms, and pyramids.

\subsection*{Tetrahedron}

\

\begin{theorem}\label{thm:3d-t}
Let $K$ be the reference tetrahedron with coordinates $(x,y,z)$.
Then,
the following two sequences on $K$ are  exact for $k\ge 0$,

\begin{tabular}{l c c c c c c c}
$\satet:$ &
$ \pol_{k+3}
$&$\overset{\grads}{\longrightarrow}$ &
$ \bpol_{k+2}
$&
$\overset{\curls}{\longrightarrow}$ &
$\bld\pol_{k+1}
$&
$\overset{\divs}{\longrightarrow}$ &
$\pol_{k}
$,\\
\\
$\sbtet:$ &$\pol_{k+1}
$&
$\overset{\grads}{\longrightarrow}$ &
$\bld\pol_{k}
\oplus\bld x
\times\bld{\widetilde\pol}_{k}
$
&
$\overset{\curls}{\longrightarrow}$ &
$\bpol_{k}
\oplus \bld x\,\widetilde\pol_k
$&
$\overset{\curls}{\longrightarrow}$ &
$\pol_{k}
.$
\end{tabular}\\
Here $\bld x=(x,y,z)^t$ and $\pol_k$ denotes $\pol_k(x,y,z)$.

Moreover, the sequence of traces  on a face $\fc$ for $\satet$ is 
the sequence $\satrip(\fc)$, and  that for $\sbtet$ is
the sequence $\sbtri(\fc)$.
\end{theorem}

\

Both of these sequences are well-known.
The first sequence $\satet$ is mainly due to N\'ed\'elec \cite{Nedelec86} since 
its $H(\mathrm{curl})$ and $H(\mathrm{div})$ spaces are N\'ed\'elec's edge and face spaces 
of second kind of degree $k+2$ and $k+1$, respectively. 
Its $H^1$ and $L^2$ spaces are the Lagrange polynomial space of degree $k+3$, and 
the discontinuous polynomial space of degree $k$.
The second sequence $\sbtet$ is mainly due to N\'ed\'elec \cite{Nedelec80} since 
its $H(\mathrm{curl})$ and $H(\mathrm{div})$ spaces are N\'ed\'elec's edge and face spaces 
of first kind of degree $k$, respectively. 
Its $H^1$ and $L^2$ spaces are the Lagrange polynomial space of degree $k+1$, and the
discontinuous polynomial space of degree $k$, respectively.

\subsection*{Hexahedron}

\

\begin{theorem}\label{thm:3d-hex}
Let $K$ be the reference hexahedron with coordinates $(x,y,z)$. Then,
the following four sequences are exact for $k\ge 0$:

\resizebox{0.95\columnwidth}{!}{
\begin{tabular}{l c c c c c c c}
$\sacub:$ &
$ \pol_{k+3}
$&$\overset{\grads}{\longrightarrow}$ &
$ \bpol_{k+2}
$&
$\overset{\curls}{\longrightarrow}$ &
$\bld\pol_{k+1}
$&
$\overset{\divs}{\longrightarrow}$ &
$\pol_{k}
$,\\
& $\oplus\delta H^{3,I}_{k+3}$ &
& $\oplus\grads\delta H^{3,I}_{k+3}$&
&$\oplus\curls\delta E^{3,I}_{k+2}$\\
&  &
& $\oplus\delta E^{3,I}_{k+2}$\\
$\sbcub:$ &
$ \pol_{k+1}
$&$\overset{\grads}{\longrightarrow}$ &
$ \bpol_{k}\oplus \bld x\times \widetilde\bpol_k
$&
$\overset{\curls}{\longrightarrow}$ &
$\bld\pol_{k}\oplus \bld x\,\widetilde\pol_k
$&
$\overset{\divs}{\longrightarrow}$ &
$\pol_{k}
$,\\
& $\oplus\delta H^{3,I}_{k+1}$ &
& $\oplus\grads\delta H^{3,I}_{k+1}$&
&$\oplus\curls\delta E^{3,I}_{k+1}$\\
&  &
& $\oplus\delta E^{3,I}_{k+1}$
\vspace{.3cm}\\
\hline
\\
$\sccub:$: &
$ \qol_{k}\oplus \left\{
\begin{tabular}{c}
 $x^{k+1},$\\$y^{k+1},$
 \\
 $z^{k+1}$
\end{tabular}
\right\}
$&$\overset{\grads}{\longrightarrow}$ &
$ \bqol_{k}\oplus \bld x\times 
\left\{
\begin{tabular}{c}
$ y^kz^k\grads x$, \\
$z^kx^k\grads y,$\\
$x^ky^k\grads z$
\end{tabular}
\right\}
$&
$\overset{\curls}{\longrightarrow}$ &
$\bld\qol_{k}\oplus \bld x\,\{x^ky^kz^k\}
$&
$\overset{\divs}{\longrightarrow}$ &
$\qol_{k}
$,\\
& $\oplus\delta H^{3,II}_{k+1}$ &
& $\oplus\grads\delta H^{3,II}_{k+1}$&
&$\oplus\curls\delta E^{3,II}_{k+1}$\\
&  &
& $\oplus\delta E^{3,II}_{k+1}$\\
$\sdcub:$ &
$ \qol_{k+1}$&
$\overset{\grads}{\longrightarrow}$ &
$ 
\left(
\begin{tabular}{c}
$ \pol_{k,k+1,k+1}$ \\
$\pol_{k+1,k,k+1}$\\
$\pol_{k+1,k+1,k}$
\end{tabular}
\right)
$&
$\overset{\curls}{\longrightarrow}$ &
$\left(
\begin{tabular}{c}
$ \pol_{k+1,k,k}$ \\
$\pol_{k,k+1,k}$\\
$\pol_{k,k,k+1}$
\end{tabular}
\right)
$&
$\overset{\divs}{\longrightarrow}$ &
$\qol_{k}
$.\\
\end{tabular}
}\\
Here the additional spaces $\delta H^{3,I}_{k},\delta H^{3,II}_{k}$, 
and $\delta E^{3,I}_{k},\delta E^{3,II}_{k}$ 
for $k\ge 1$
takes the following forms: 
\begin{subequations}
\begin{alignat*}{2}
 \delta H^{3,I}_{k} :=&\;\mathrm{span} \left\{
 \begin{tabular}{c}
   $x\,y\,z^{k},\;y\,z\,x^k,\;z\,x\,y^k,$\\
   $x\,\widetilde\pol_{k}(y,z),\;
   y\,\widetilde\pol_{k}(z,x),\;
   z\,\widetilde\pol_{k}(x,y)$\\
 \end{tabular}
 \right\},\\
 \delta H^{3,II}_{k} :=&\;\mathrm{span} \left\{
 \begin{tabular}{c}
   $x\,y\,z^{k},\;y\,z\,x^k,\;z\,x\,y^k,$\\
   $x\,y^k,\; y\,z^k,\; z\,x^k$,\\
   $x\,z^k,\; y\,x^k,\; z\,y^k$\\
 \end{tabular}
 \right\},\\ 
 \end{alignat*}
\begin{alignat*}{2}
 \delta E^{3,I}_{k} :=&\;\mathrm{span} \left\{
 \begin{tabular}{c}
   $x\,\widetilde\pol_{k-1}(y,z)(y\,\grads z-z\,\grads y)$,\\
$y\,\widetilde\pol_{k-1}(z,x)(z\,\grads x-x\,\grads z)$,\\
$z\,\widetilde\pol_{k-1}(x,y)(x\,\grads y-y\,\grads x)$
 \end{tabular}
 \right\},\\
 \delta E^{3,II}_{k} :=&\;
 \mathrm{span} \left\{
 \begin{tabular}{c}
   $x(y^k\,\grads z-z^k\,\grads y)$,
$y(z^k\,\grads x-x^k\,\grads z)$,\\
$z(x^k\,\grads y-y^k\,\grads x)$,
  $x\,y^{k-1}\,z^{k-1}(y\,\grads z-z\,\grads y)$,\\
$y\,z^{k-1}\,x^{k-1}(z\,\grads x-x\,\grads z)$,
$z\,x^{k-1}\,y^{k-1}(x\,\grads y-y\,\grads x)$
 \end{tabular}
 \right\}. 
\end{alignat*}
\end{subequations}

Moreover, 
the sequence of traces on a square face $\fc$
for $\sacub$ is the sequence $\sasqrp(\fc)$,
and  that for $\sicub$ is the sequence $\sisqr(\fc)$ 
for $i\in\{2,3,4\}$.
\end{theorem}

\

Note that for $k=0$, the last three sequences are exactly the same.
Here the second sequence is new; the other three are known.
The first sequence $\sacub$
is the 
serendipity sequence of Arnold and Awanou \cite{ArnoldAwanou14}.
The second sequence is a slight variation of the first one. It displays  a new family of 
$H(\mathrm{curl})$
and $H(\mathrm{div})$ spaces.
The  third sequence $\sccub$ is the TNT sequence of Cockburn and Qiu 
 \cite{CockburnQiu14} 
(with a slight variation in the space representation).
The last sequence is mainly due to N\'ed\'elec \cite{Nedelec80}
since its $H(\mathrm{curl})$ and $H(\mathrm{div})$ spaces are N\'ed\'elec's edge and face spaces 
of first kind of degree $k$ on the 
cube, respectively.  Its $H^1$ and $L^2$ spaces are the tensor-product Lagrange polynomial space of degree $k+1$, and 
the discontinuous tensor-product polynomial space of degree $k$, respectively.

\subsection*{Prism}

\

\begin{theorem}\label{thm:3d-prism}
Let $K$ be the reference prism with coordinates $(x,y,z)$. Then,
the following four sequences are exact for $k\ge 0$:

\resizebox{0.95\columnwidth}{!}{
\begin{tabular}{l c c c c c c c}
$\sapri:$ &
$ \pol_{k+3}
$&$\overset{\grads}{\longrightarrow}$ &
$ \bpol_{k+2}
$&
$\overset{\curls}{\longrightarrow}$ &
$\bld\pol_{k+1}
$&
$\overset{\divs}{\longrightarrow}$ &
$\pol_{k}
$,\\
& $\oplus\delta H^{3,III}_{k+3}$ &
& $\oplus\grads\delta H^{3,III}_{k+3}$&
&$\oplus\curls\delta E^{3,III}_{k+2}$\\
&  &
& $\oplus\delta E^{3,III}_{k+2}$\\
$\sbpri:$ &
$ \pol_{k+1}
$&$\overset{\grads}{\longrightarrow}$ &
$ \bpol_{k}\oplus \bld x\times \widetilde\bpol_k
$&
$\overset{\curls}{\longrightarrow}$ &
$\bld\pol_{k}\oplus \bld x\,\widetilde\pol_k
$&
$\overset{\divs}{\longrightarrow}$ &
$\pol_{k}
$,\\
& $\oplus\delta H^{3,III}_{k+1}$ &
& $\oplus\grads\delta H^{3,III}_{k+1}$&
&$\oplus\curls\delta E^{3,III}_{k+1}$\\
&  &
& $\oplus\delta E^{3,III}_{k+1}$
\vspace{.3cm}\\
\hline
\\
 $\scpri:$ &
  $ \pol_{k|k}
 $&
 $\overset{\grads}{\longrightarrow}$ &
 $ \bpol_{k|k}
$
&
$\overset{\curls}{\longrightarrow}$ &
$\bld\pol_{k|k}$&
$\overset{\divs}{\longrightarrow}$ &
$\pol_{k|k}
$,\\
&\!\!\!\!\!\!\!\!
 $\oplus\widetilde\pol_{k+1}(x,y)\oplus
\{ z^{k+1}\}
 $
 &
 &
 $\oplus 
 \left(\!\!\!\!
 \begin{tabular}{c}
$  \left(\!\!\!\!
  \begin{tabular}{c}
 $y$\\
 $-x$
 \end{tabular}
\!\!\!\! \right)\widetilde\pol_{k}(x,y)$\\
$\widetilde\pol_{k+1}(x,y)\,z^k$
\end{tabular}
\!\!\!\! \right)
$
&
&
$\oplus 
 \left(\!\!\!\!
 \begin{tabular}{c}
$0$\\
$0$\\
$\widetilde\pol_{k}(x,y)\,z^{k+1}$
\end{tabular}
\!\!\!\! \right)
$&\\
& $\oplus\delta H^{3,III}_{k+1}$ &
& $\oplus\grads\delta H^{3,III}_{k+1}$&
&$\oplus\curls\delta E^{3,III}_{k+1}$\\
&  &
& $\oplus\delta E^{3,III}_{k+1}$\\
$\sdpri:$ &
$ \pol_{k+1|k+1}$&
$\overset{\grads}{\longrightarrow}$ 
&
$ 
 \left(\!\!\!\!
 \begin{tabular}{c}
$
\bld{N}_{k}(x,y)\otimes\pol_{k+1}(z)$\\	
$\pol_{k+1|k}$
\end{tabular}
\!\!\!\! \right)
$&
$\overset{\curls}{\longrightarrow}$ &
$\left(\!\!\!\!
 \begin{tabular}{c}
$ \bld{{RT}}_{k}(x,y)\otimes\pol_{k}(z)$\\	
$\pol_{k|k+1}$
\end{tabular}
\!\!\!\! \right)
$&
$\overset{\divs}{\longrightarrow}$ &
$\pol_{k|k}
$.\\
\end{tabular}
}\\
Here, we have 
\[
 \bld N_k(x,y): = \bpol_k(x,y)\oplus \left(\hspace{-.2cm}\begin{tabular}{c}
                                      $-y$\\
                                      $x$
                                     \end{tabular}\hspace{-.2cm}\right)\widetilde{\pol}_k(x,y),
                                     \quad
 \bld {RT}_k(x,y): = \bpol_k(x,y)\oplus \left(\hspace{-.2cm}\begin{tabular}{c}
                                      $x$\\
                                      $y$
                                     \end{tabular}\hspace{-.2cm}\right)\widetilde{\pol}_k(x,y),                                     
\]
and
the additional spaces $\delta H^{3,III}_{k}$
and $\delta E^{3,III}_{k}$ 
for $k\ge 1$
takes the following forms: 
\begin{alignat*}{2}
 \delta H^{3,III}_{k} :=&\;\mathrm{span} \left\{
 \begin{tabular}{c}
   $z^{k}\,\widetilde\pol_1(x,y),\;z\,\widetilde\pol_k(x,y)$\\
 \end{tabular}
 \right\},\\
 \delta E^{3,III}_{k} :=&\;\mathrm{span} \left\{
 \begin{tabular}{c}
   $z^{k}(x\,\grads y-y\,\grads x),\;
   z\,\widetilde\pol_{k-1}(x,y)(x\,\grads y-y\,\grads x)$
 \end{tabular}
 \right\}.
\end{alignat*}

Moreover, the sequence of traces for these four sequences are
\begin{alignat*}{3}
 \tr^\fc\left(
 \sapri\right) =&\; \left\{
 \begin{tabular}{l l}
 $\satrip(\fc)$ & triangle $\fc$,
 \vspace{.2cm}\\
 $\sasqrp(\fc)$ & square $\fc$,
 \end{tabular}
 \right.
 &&\quad
  \tr^\fc\left(
 \sbpri\right) =&&\; \left\{
 \begin{tabular}{l l}
 $\sbtri(\fc)$ & triangle $\fc$,
 \vspace{.2cm}\\
 $\sbsqr(\fc)$ & square $\fc$,
 \end{tabular}
 \right.
 \\
 \tr^\fc\left(
 \scpri\right) =&\; \left\{
 \begin{tabular}{l l}
 $\sbtri(\fc)$ & triangle $\fc$,
 \vspace{.2cm}\\
 $\widetilde{\scsqr}(\fc)$ & square $\fc$,
 \end{tabular}
 \right.
&& \quad
 \tr^\fc\left(
 \sapri\right) =&&\; \left\{
 \begin{tabular}{l l}
 $\sbtri(\fc)$ & triangle $\fc$,
 \vspace{.2cm}\\
 $\sdsqr(\fc)$ & square $\fc$,
 \end{tabular}
 \right.
\end{alignat*}
where the sequence  
$
 \widetilde{\scsqr}(\fc)
 $
 on the square face with coordinates $(\xi,z)$
is a slight 
modification of the sequence  ${\scsqr}(\fc)$, with the same $H^1$ and $L^2$ spaces,
but a different $H(\mathrm{curl})$ space, namely,
$
 \bqol_k(\xi,z)\oplus \xi^{k+1}z^k\grads z\oplus
 \grads \{\xi\,z^{k+1},z\,\xi^{k+1}\}.
$
\end{theorem}

\

Note that for $k=0$, the sequences $\sbpri$ and $\sdpri$
are exactly the same; they are slightly different from the sequence $\scpri$.
Note also that the $H(\mathrm{curl})$ space
for $
 \widetilde{\scsqr}(\fc)
 $
is not invariant under the coordinate permutation $(\xi, z)\rightarrow (z,\xi)$,
where $(\xi,z)$ are the coordinates on $\fc$.
We did not find a prismatic sequence with trace space 
on a square face $\fc$ exactly equal
to $\scsqr(\fc)$. This brings about a small complication for constructing
the $H(\mathrm{curl})$-conforming finite element spaces on a hybrid mesh, 
see the discussion on commuting exact sequences on the whole polyhedral mesh
below. 

Finally, note that, while the fourth sequence is known,
see \cite{FuentesKeithDemkowiczNagaraj15},
the other three are new.
These four sequences can be considered as the extensions of the related sequences 
on the reference cube to the
reference prism. 

\subsection*{Pyramid}

\

\begin{theorem}\label{thm:3d-pyramid}
Let $K$ be the reference pyramid with coordinates $(x,y,z)$. Then,
the following four sequences are exact for $k\ge 0$:

\

\resizebox{0.95\columnwidth}{!}{
\begin{tabular}{l c c c c c c c}
$\sapyr:$ &
$ \pol_{k+3}
$&$\overset{\grads}{\longrightarrow}$ &
$ \bpol_{k+2}
$&
$\overset{\curls}{\longrightarrow}$ &
$\bld\pol_{k+1}
$&
$\overset{\divs}{\longrightarrow}$ &
$\pol_{k}
$,\\
& $\oplus\delta H^{3,IV}_{k+3}$ &
& $\oplus\grads\delta H^{3,IV}_{k+3}$&
&$\oplus\curls\delta E^{3,IV}_{k+2}$\\
&  &
& $\oplus\delta E^{3,IV}_{k+2}$\\
$\sbpyr:$ &
$ \pol_{k+1}
$&$\overset{\grads}{\longrightarrow}$ &
$ \bpol_{k}\oplus \bld x\times \widetilde\bpol_k
$&
$\overset{\curls}{\longrightarrow}$ &
$\bld\pol_{k}\oplus \bld x\,\widetilde\pol_k
$&
$\overset{\divs}{\longrightarrow}$ &
$\pol_{k}
$,\\
& $\oplus\delta H^{3,IV}_{k+1}$ &
& $\oplus\grads\delta H^{3,IV}_{k+1}$&
&$\oplus\curls\delta E^{3,IV}_{k+1}$\\
&  &
& $\oplus\delta E^{3,IV}_{k+1}$
\vspace{.3cm}\\
\hline
\\
$\scpyr:$ &
$ \pol_{k+1}
$&$\overset{\grads}{\longrightarrow}$ &
$ \bpol_{k}\oplus \bld x\times \widetilde\bpol_k
$&
$\overset{\curls}{\longrightarrow}$ &
$\bld\pol_{k}\oplus \bld x\,\widetilde\pol_k
$&
$\overset{\divs}{\longrightarrow}$ &
$\pol_{k}
$,\\
& $\oplus\delta H^{3,V}_{k+1}$ &
& $\oplus\grads\delta H^{3,V}_{k+1}$&
&$\oplus\curls\delta E^{3,V}_{k+1}$\\
&  &
& $\oplus\delta E^{3,V}_{k+1}$\\
$\sdpyr:$ &
$ \pol_{k+1}
$&$\overset{\grads}{\longrightarrow}$ &
$ \bpol_{k}\oplus \bld x\times \widetilde\bpol_k
$&
$\overset{\curls}{\longrightarrow}$ &
$\bld\pol_{k}\oplus \bld x\,\widetilde\pol_k
$&
$\overset{\divs}{\longrightarrow}$ &
$\pol_{k}
$,\\
& $\oplus\delta H^{3,VI}_{k+1}$ &
& $\oplus\grads\delta H^{3,VI}_{k+1}$&
&$\oplus\curls\delta E^{3,V}_{k+1}$\\
&  &
& $\oplus\delta E^{3,V}_{k+1}$\\
\end{tabular}
}\\
Here the additional spaces $\delta H^{3,IV}_{k},\delta H^{3,V}_{k},
\delta H^{3,VI}_{k}$
and $\delta E^{3,IV}_{k},\delta E^{3,V}_{k},
\delta E^{3,VI}_{k}$ 
for $k\ge 1$
takes the following forms: 
\begin{alignat*}{2}
 \delta H^{3,IV}_{k} :=&\;\mathrm{span} \left\{
 \begin{tabular}{c}
   $\frac{x\,y}{1-z}z^{k-1},
   \frac{x\,y\,z}{1-z}\widetilde\pol_{k-2}(x,z),
   \frac{x\,y\,z}{1-z}\widetilde\pol_{k-2}(y,z)$\\
 \end{tabular}
 \right\},\\
 \delta H^{3,V}_{k} :=&\;\delta H^{3,IV}_{k}\oplus
 \mathrm{span} \left\{
\hspace{-.2cm} \begin{tabular}{c}
   $\frac{x^\alpha\,y^\beta}{(1-z)^{\min\{\alpha,\beta\}}}:\!\!\!\!
\begin{tabular}{c}
$ \alpha=1, \beta=k$ or 
$ \alpha=k, \beta=1$\\
or $\alpha \le k-1,\; \beta \le k-1,\; k+1\le \alpha+\beta$
\end{tabular}
$\\
 \end{tabular}
\hspace{-.3cm} \right\}
 ,\\
  \delta H^{3,VI}_{k} :=&\;\delta H^{3,IV}_{k}\oplus
 \mathrm{span} \left\{
 \begin{tabular}{c}
   $\frac{x^\alpha\,y^\beta}{(1-z)^{\min\{\alpha,\beta\}}}:
   \;\;
\begin{tabular}{c}
$ \alpha \le k,\;  \beta \le k,\; k+1\le \alpha+\beta$
\end{tabular}
$\\
 \end{tabular}
 \right\},\\
 \delta E^{3,IV}_{k} :=&\;\mathrm{span} \left\{
 \begin{tabular}{c}
   $\frac{x\,y^k}{1-z}\grads z,\frac{y\,x^k}{1-z}\grads z,
   \frac{x\,y\, z}{1-z}\grads x$
 \end{tabular}
 \right\},\\
  \delta E^{3,V}_{k} :=&\;
   \delta E^{3,IV}_{k}\oplus\mathrm{span} 
   \left\{
  \begin{tabular}{c}
    $\frac{x^\alpha\,y^{\beta+1}}{(1-z)^{\min\{\alpha,\beta\}}}\grads x:
    \;\;
 \begin{tabular}{c}
 $\alpha \le k,\;
 \beta \le k,\; k+1\le \alpha+\beta$
 \end{tabular}
$
 \end{tabular}
 \right\}.
\end{alignat*}

Moreover, the sequence of traces for these four sequences are
\begin{align*}
 \tr^\fc\left(
 \sapyr\right) =&\; \left\{
 \begin{tabular}{l l}
 $\satrip(\fc)$ & if $\fc$ is a triangle,
 \vspace{.2cm}\\
 $\sasqrp(\fc)$ & if $\fc$ is the base square,
 \end{tabular}
 \right.\\
 \tr^\fc\left(
 \sipyr\right) =&\; \left\{
 \begin{tabular}{l l}
 $\sbtri(\fc)$ & if $\fc$ is a triangle,
 \vspace{.2cm}\\
 $\sisqr(\fc)$ & if $\fc$ is the base square,
 \end{tabular}
 \right. & \text{ for }i\in\{2,3,4\}.
\end{align*}
\end{theorem}

\

All the four sequences are new and contain certain rational functions.
They can be considered as the extensions of the related sequences 
on the reference cube to the 
reference pyramid. 
Note that, unlike all the previous sequences,  the function spaces here includes certain 
rational functions.

{
We remark that similar {\em serendipity-type} pyramidal $H^1$-conforming space as that in $\sapyr$
(with the same space dimension) 
was recently introduced in \cite{Gillette15}. 
The serendipity space in \cite{Gillette15}, containing $\pol_k(K)$,
was obtained by mapping certain rational functions from the {\em infinite pyramid} to the reference pyramid, extending similar results in \cite{NigamPhillips12,NigamPhillips12b} with a significant dimension reduction.
A set of degrees of freedom, similar as those in \cite{ArnoldAwanou14} on cubes, was also identified. 
We believe the $H^1$-conforming space in \cite{Gillette15} (obtained from the infinite pyramid) and that in $\sapyr$ (obtained directly from the reference pyramid) should be closely related since both enrich same amount of rational functions (whose precise definition seems to be different) to $\pol_k(K)$ to achieve conformity. 
}

Let us point out that our fourth sequence $\sdpyr$ 
is significantly smaller than the related pyramidal sequence presented in \cite{FuentesKeithDemkowiczNagaraj15},
originally obtained in \cite{NigamPhillips12}.

\subsection*{Whole mesh} 
 Using the sequences in Theorem \ref{thm:3d-t} to 
 Theorem \ref{thm:3d-pyramid}, 
 we can readily obtain {\it mapped} 
 commuting exact sequences on a hybrid mesh, $\Omega_h:=\{K\}$, of 
a polyhedral domain $\Omega$, where each physical element $K$ is 
an affine mapping of any of the 
four  reference polyhedra.
We obtain {\it four} family of {\it mapped } 
sequences with the following form,  
\[
\mathrm{S_{3d}}(\Omega_h):\quad  H_3(\Omega_h)
\overset{\grads}{\longrightarrow} E_3(\Omega_h)
 \overset{\curls}{\longrightarrow} V_3(\Omega_h)
 \overset{\divs}{\longrightarrow} W_3(\Omega_h),
\]
where $H_3(\Omega_h)\times E_3(\Omega_h)
 \times V_3(\Omega_h)\times W_3(\Omega_h)\subset 
 H^1(\Omega)\times H(\mathrm{curl},\Omega)
 \times H(\mathrm{div},\Omega)\times L^2(\Omega)$.

The restriction on an element $K$ 
of 
the first family of sequences 
is $\satet(K)$ if $K$ is a tetrahedron, 
$\sacub(K)$ if $K$ is a parallelepiped,
$\sapri(K)$ if $K$ is a parallel prism, 
and
$\sapyr(K)$ if $K$ is a pyramid with a parallelogram base.

The restriction on an element $K$ 
of 
the second family of sequences 
is $\sbtet(K)$ if $K$ is a tetrahedron, 
$\sbcub(K)$ if $K$ is a parallelepiped,
$\sbpri(K)$ if $K$ is a parallel prism, 
and
$\sbpyr(K)$ if $K$ is a pyramid with a parallelogram base.

The restriction on an element $K$ 
of 
the third family of sequences, 
defined on a hybrid mesh without prisms,
is 
$\sbtet(K)$ if $K$ is a tetrahedron, 
$\sccub(K)$ if $K$ is a parallelepiped,
and
$\scpyr(K)$ if $K$ is a pyramid with a parallelogram base. 

The restriction on an element $K$ 
of  the fourth family of sequences 
is $\sbtet(K)$ if $K$ is a tetrahedron, 
$\sdcub(K)$ if $K$ is a parallelepiped,
$\sdpri(K)$ if $K$ is a parallel prism, 
and
$\sdpyr(K)$ if $K$ is a pyramid with a parallelogram base.
Note that this family of sequences is a modification 
of the one considered in \cite{FuentesKeithDemkowiczNagaraj15}
on a hybrid mesh with a smaller pyramidal sequence.

The reason we are able to do so is due to 
the trace compatibilities in Theorem \ref{thm:3d-t} to
Theorem \ref{thm:3d-pyramid}. 
In particular, we mention that we 
exclude prisms in the hybrid mesh of
the third family of sequences mainly due to 
the  incompatibility of the $H(\mathrm{curl})$ trace of $\scpri$ and that of 
$\sccub$, which differ by a single function.

\section{Proofs of the main results in Section \ref{sec:mainresults}}
\label{sec:proof}

\subsection{Proof of Theorem \ref{thm:1d-general}}
\begin{proof}
Suppose the segment $K$ has coordinate $x$.
By the exactness of the sequence $S_1^g(K)$, we have 
$\bR = \kerl_{\grads} H_1^g(K)$ (and hence 
$1 \in H_1^g(K)$), and 
$\pol_0(K)\subset W_1^g(K) = \grads H_1(K)$ (and hence 
$x\in H_1^g(K)$). This implies that 
$\pol_1(K)\subset H_1^g(K)$, and hence 
$\dim \tr_H(H_1^g(K)) = \dim \tr_H(\pol_1(K)) = 2$.
So $S_1^g(K)$ is a compatible exact sequence.
\end{proof}


\subsection{Proof of Theorem \ref{thm:2d-general}}
\begin{proof}
To simplify the notation, 
we set $H_2(K):= H_2^g(K)\oplus \delta H_2^g(K)$,
$E_2(K):= E_2^g(K)\oplus \grads\delta H_2^g(K)$, and 
$W_2(K):= W_2^g(K)$.

The exactness of the sequence $S_2(K)$ follows directly from
the exactness of the given sequence
\[
  H_2^g(K)\overset{\grads}{\longrightarrow}E_2^g(K)
 \overset{\curls}{\longrightarrow}W_2^g(K).
\]

It is easy to show that for any 
exact sequence $H_2(K)\rightarrow E_2(K)\rightarrow W_2(K)$, we have
\begin{subequations}
\label{2d-claim} 
 \begin{align}
\label{2d-claim-0} 
\dim H_2(K) -
 \dim E_2(K) +\dim W_2(K)\; = 1,\\
\label{2d-claim-1} 
 \grads \oo{H_2} = \; \{v\in \oo{E_2}:\;\;\curls v=0\},\hspace{1cm}\\
 \label{2d-claim-2} 
 \curls \oo{E_2}(K) \subset\; \oo{W_2}(K).\hspace{1.7cm} 
\end{align}
 \end{subequations}

In view of Definition \ref{2d-cex} on
a compatible exact sequence, we need to prove  that
the two dimension count identities of
Property (ii) are satisfied, and that the sequence of traces
on each edge is a compatible exact sequence.

We have
\begin{alignat}{2}
\label{oc-0}
 \dim \tr_H H_2(K) = &\;
 \dim H_2(K) - \dim \oo{H_2}(K)\\
=&\; \dim H_2^g(K) +
 \dim \delta H_2^g(K)- \dim \oo{H_2}(K)
 \nonumber\\
=&\; \dim H_2^g(K) +
 \dim \delta H_2^g(K)
- \dim \oo{H_2^g}(K) 
\nonumber\\
=&\; \sum_{\eg\in\EG(K)} (\dim \oo{H_1}(\eg) +1),
\nonumber
\end{alignat}
where the second equality is due to property 
(ii) of $\delta H_2^g(K)$, the third one is due to 
property (iii), and the last one is due to property (iv). 
Now, by property (i) of the
Definition \ref{s1-ad}
of an $S_1(\dK)$ admissible
exact sequence, and by property (i) of $\delta H_2^g(K)$, 
we have 
\begin{align*}
\tr_H^\eg H_2(K) \subset H_1(\eg)\quad \forall \eg\in\EG(K). 
\end{align*}
On the other hand, since 
\[
\left. \tr_HH_2(K)\right|_\eg
 \subset \tr_H^\eg H_2(K)
 \quad\forall \eg\in\EG(K),
 \]
we have 
\begin{align*}
 \dim \tr_HH_2(K) 
 \le&\; \sum_{\eg\in\EG(K)} \dim
 \oo{\overline{\tr_H^\eg H_2}}(K)
 +\sum_{\vt\in\VT(K)} 1
 \le 
 \sum_{\eg\in\EG(K)} (\dim
 {\oo{H_1}(\eg)}+1).
\end{align*}
By \eqref{oc-0}, we have the above inequalities are indeed
equalities, and 
\begin{align}
\label{oc-1}
\tr_H^\eg H_2(K) = H_1(\eg)\quad \forall \eg\in\EG(K). 
\end{align}
This proves the first dimension count identity of 
Property (ii) for a compatible exact sequence.

Next, let us prove the second dimension count identity.
We have that 
\begin{alignat}{2}
\label{oc-2}
 \dim \tr_E E_2(K) = &\;
 \dim E_2(K) - \dim \oo{E_2}(K)\\
=&\; \dim H_2(K)+
\dim W_2(K) -1 \nonumber\\
&\;- 
\dim \curls\oo{E_2}(K)
-\dim \{v\in\oo{E_2}(K):\;
\curls v =0
\}
 \nonumber\\
=&\; \dim H_2(K) +
 \dim \oo{W_2}(K)
- \dim \curls\oo{E_2}(K) 
-\dim \oo{H_2}(K)
\nonumber\\
=&\; \sum_{\eg\in\EG(K)} (\dim \oo{H_1}(\eg) +1)
+ \dim \oo{W_2}(K)
- \dim \curls\oo{E_2}(K) 
\nonumber\\
\ge&\;\sum_{\eg\in\EG(K)} (\dim \oo{H_1}(\eg) +1)\nonumber\\
=&\;\sum_{\eg\in\EG(K)} \dim W_1(\eg) \nonumber,
\end{alignat}
where the second equality is due to \eqref{2d-claim-0} of an exact sequence,
the third equality is due to 
\eqref{2d-claim-1} and Property (ii) of Definition 
\ref{s1-ad}, the fourth equality is due to \eqref{oc-0},
the fifth inequality is due to \eqref{2d-claim-2}, and 
the last equality is due to the exactness of the sequences
$S_1(\eg)$.
Now, by property (i) of 
Definition \ref{s1-ad}
for an $S_1(\dK)$ admissible
exact sequence, 
we have 
\begin{align*}
\tr_E^\eg E_2^g(K) \subset W_1(\eg)\quad \forall \eg\in\EG(K), 
\end{align*}
and by  property (i) of $\delta H_2^g(K)$, we have 
\begin{align*}
\tr_E^\eg \grads \delta H_2^g(K)
=
\grads\tr_H^\eg \delta H_2^g(K)
\subset \grads H_1(\eg) = W_1(\eg)\quad \forall \eg\in\EG(K).
\end{align*}
Hence, 
\begin{align*}
\tr_E^\eg E_2(K) \subset W_1(\eg)\quad \forall \eg\in\EG(K). 
\end{align*}

On the other hand, since 
\[
\left. \tr_EE_2(K)\right|_\eg
 \subset \tr_E^\eg E_2(K)
 \quad\forall \eg\in\EG(K),
 \]
we have 
\begin{align*}
 \dim \tr_EE_2(K) 
 \le&\; \sum_{\eg\in\EG(K)} \dim
 {{\tr_E^\eg E_2}}(K)
\le
 \sum_{\eg\in\EG(K)} \dim
 W_1(\eg).
\end{align*}
By \eqref{oc-2}, we have the above inequalities are indeed 
equalities.
Hence, 
\begin{align}
\label{oc-5}
\tr_E^\eg E_2(K) = W_1(\eg)
 \quad\forall \eg\in\EG(K). 
\end{align}
This completes the proof of the second dimension count
identity of Property (ii) for a compatible exact sequence in Definition \ref{2d-cex}.
The equalities \eqref{oc-1} and \eqref{oc-5} 
ensure that the sequence of traces 
$\tr^\eg(S_2(K))=S_1(\eg)$ is a compatible exact sequence, for all edges $\eg\in\EG(K)$. So $S_2(K)$ is a compatible exact
sequence.

Finally, invoking \cite[Corollary 3.2]{ChristiansenGillette15},
we get the minimality of the sequence $S_2(K)$ by the simple observation that  
\[
\oo{H_2}= \oo{H_2^g}, \quad
\oo{E_2}= \oo{E_2^g},\quad
\oo{W_2}= \oo{W_2^g}.
\]

The proof of Theorem  \ref{thm:2d-general} is now complete.
\end{proof}

\subsection{Proof of Theorem \ref{thm:3d-general}}
\begin{proof}
To simplify the notation, 
we set $H_3(K):= H_3^g(K)\oplus \delta H_3^g(K)$,
$E_3(K):= E_3^g(K)\oplus \grads\delta H_3^g(K)\oplus 
\delta E_3^g(K)$, 
$V_3(K):= V_3^g(K)\oplus \curls \delta E_3^g(K)$, and 
$W_3(K):=W_3^g(K)$.

The exactness of the sequence $S_3(K)$ comes directly from
the exactness of the given sequence
\[
 H_3^g(K)\overset{\grads}{\longrightarrow}E_3^g(K)
 \overset{\curls}{\longrightarrow}V_3^g(K)
 \overset{\divs}{\longrightarrow}W_3^g(K)
\]

In view of Definition \ref{3d-cex} of a compatible
exact sequence, we need to prove that
the three dimension count identities of
Property (ii) are satisfied, and that the
sequence of traces on any face $\fc\in\FC(K)$
is a compatible exact sequence.

The proof of the first dimension count identity for 
the $H^1$-trace space $\tr_HH_3(K)$
and, consequently,
that of the identity
\begin{align}
 \label{3d-h1}
 \tr_H^\fc H_3(K) = H_2(K) \quad\forall \fc\in\FC(K)
\end{align}
are omitted because they are 
very similar to those of the two-dimensional case; 
see details in the proof of Theorem \ref{thm:2d-general}.

Now, let us  prove the third dimension count identity 
for the $H(\mathrm{div})$-trace space $\tr_VV_3(K)$
of  Property (ii) in Definition \ref{3d-cex}.

By properties (ii), (iii), and 
(iv)
of $\delta E_3^g(K)$, we have
\begin{alignat}{2}
 \label{oh-v}
 \dim\tr_V V_3(K) = &\; \dim V_3(K) -\dim \oo{V_3}(K)
\\
=&\;\dim V_3^g(K)+ \dim \delta E_3^g(K) 
 -\dim \curls \oo{V_3}(K)\nonumber\\
&\;-\dim \{v\in \oo{V_3}(K):\curls v=0\} 
\nonumber
 \\
=&\; 
\sum_{\fc\in\FC(K)}\dim {W_2}(\fc) 
+\dim \oo{W_3^g}(K)
-\dim \curls \oo{V_3}(K)
\nonumber
\\
\ge &\;
\sum_{\fc\in\FC(K)}\dim {W_2}(\fc),
\nonumber
\end{alignat}
where in the last inequality,  we used the fact that 
\[
 \divs \oo{V_3}(K)\subset \oo{W_3}(K) = \oo{W_3^g}(K).
 \]

Now, by Property (i) of Definition \ref{s2-ad}
of an $S_2(\dK)$-admissible exact sequence we have 
$\tr_V^\fc V_3^g(K)\subset W_2(\fc)$, and
by 
property $\mathrm{(i)}$ of 
$\delta E_3^g(K)$ and 
the exactness
of the sequence $S_2(\fc)$, we have
\[\tr_V^\fc \curls \delta E_3^g(K)=
\curls\tr_E^\fc \delta E_3^g(K)\subset 
\curls E_2(\fc)=
W_2(\fc).\]
Hence, 
\[
 \tr_V^\fc V_3(K)\subset W_2(\fc)
 \quad \forall \fc\in\FC(K).
\]
This implies that 
\begin{align*}
 \dim \tr_VV_3(K)\le &\;
 \sum_{\fc\in\FC(K)}
 \tr_V^\fc V_3(K)
 \le \; \sum_{\fc\in\FC(K)} W_2(\fc).
\end{align*}
By \eqref{oh-v}, the above inequalities are indeed equalities, and
we have 
\begin{align}
\label{oh-v2}
 \tr_V^\fc V_3(K) = &\; W_2(\fc)\quad \forall \fc\in\FC(K),\\
 \label{oh-v3}
 \divs\oo{V_3}(K) = &\; \oo{W_3}(K).
\end{align}
This completes the proof of the third dimension count identity of 
Property (ii) in Definition \ref{3d-cex}.

Now, let us prove the second dimension count identity for
the $H(\mathrm{curl})$-trace space $\tr_EE_3(K)$ of Property (ii)
in Definition \ref{3d-cex}.
We use the following results of an exact sequence:
\begin{subequations}
\label{3d-edge}
\begin{align}
\label{3d-edge-0}
\dim H_3(K) -
 \dim E_3(K) +\dim V_3(K) -\dim W_3(K) =\;1,\\
\label{3d-edge-1}
\tr_E^\fc E_3(K)\subset \; E_2(\fc)\quad 
\forall \fc\in \FC(K),\hspace{1.8cm}\\
\label{3d-edge-2}
\grads\oo{H_3}(K) = \; \{v\in\oo{E_3}(K),\;\curls v = 0
\},\hspace{1.5cm}\\
\label{3d-edge-3}
\curls\oo{E_3}(K)\subset \; 
\{v\in\oo{V_3}:\;\divs v =0\}.\hspace{1.6cm}
\end{align}
\end{subequations}
Their proofs are trivial and hence omitted.
The inclusion \eqref{3d-edge-1} implies that 
\begin{align}
\label{oh-e}
 \dim \tr_EE_3(K) \le &\;\sum_{\eg\in\EG(K)}
 {{\tr_E^\eg E_3}}(K)
 +\sum_{\vt\in\VT(K)}\oo{\overline{\tr_E^\fc E_3}}(K)\\
 \le &\;\sum_{\eg\in\EG(K)}
 {W_1}(\eg)
 +\sum_{\vt\in\VT(K)}\oo{E_2}(\fc).\nonumber 
\end{align}
On the other hand, we have, by \eqref{3d-edge-0} of an exact sequence, the equalities 
\eqref{oh-v3},
\eqref{3d-edge-2}, and the inclusion
\eqref{3d-edge-3}, 
{
\begin{alignat*}{2}
 \dim \tr_EE_3(K) = &\; 
 \dim E_3(K) - \dim \oo{E_3}(K)\\
=&\; \dim E_3(K) 
- \dim \curls \oo{E_3}(K)- \dim \{v\in\oo{E_3}(K):\;\curls v =0\}
\\
=&\;
\dim H_3(K) +\dim V_3(K)-\dim W_3(K) -1\\
&\;
- \dim \curls \oo{E_3}(K)- \dim \grads 
\oo{H_3}(K) \\
= &\;
\dim \tr_HH_3(K) +\dim \tr_VV_3(K) -2\\
&\; +\dim \oo{V_3}(K)
-\dim \oo{W_3}(K)
-\dim \curls \oo{E_3}(K)\\
\ge&\;
\dim \tr_HH_3(K) +\dim \tr_VV_3(K) -2.
\end{alignat*}
}
By the first and third dimension count identities of Property (ii)
in Definition \ref{3d-cex}, we have 
the right-hand side of the above inequality is equal to
\begin{align*}
& \;\sum_{\vt\in\VT(K)}1 +
 \sum_{\eg\in \EG(K)}(\dim \oo{H_1}(\eg))
 +
 \sum_{\fc\in \FC(K)}(\dim \oo{H_2}(\fc)+\dim W_2(\fc))-2\\
= &\;
\;\sum_{\vt\in\VT(K)}1 +
 \sum_{\eg\in \EG(K)}(\dim {W_1}(\eg) -1)
 +
 \sum_{\fc\in \FC(K)}(\dim \oo{H_2}(\fc)+\dim \oo{W_2}(\fc)+1)-2\\
= & \;
\sum_{\eg\in \EG(K)}\dim {W_1}(\eg)
 +
 \sum_{\fc\in \FC(K)}\dim \oo{E_2}(\fc)- I,
\end{align*}
where 
\[
 I := \sum_{\vt\in\VT(K)} 1-\sum_{\eg\in\EG(K)} 1+
 \sum_{\fc\in\FC(K)} 1-2 = 0
\]
by  Euler's polyhedral formula.
Hence, 
\[
 \dim \tr_EE_3(K)\ge 
 \sum_{\eg\in \EG(K)}\dim {W_1}(\eg)
 +
 \sum_{\fc\in \FC(K)}\dim \oo{E_2}(\fc).
\]
And the above inequality is indeed an equality by 
\eqref{oh-e}. 
Moreover, the inclusions in \eqref{3d-edge} are also equalities: 
\begin{subequations}
\label{3d-commutingE}
\begin{align}
\label{oh-e1}
\tr_E^\fc E_3(K)=&\;
E_2(\fc)\quad\forall \fc\in \FC(K),\\
\label{oh-e2}
\curls\oo{E_3}(K)= &\; 
\{v\in\oo{V_3}:\;\divs v =0\},
\end{align} 
\end{subequations}
and this completes the proof of the second dimension count identity of 
Property (ii) in Definition \ref{3d-cex}.
The equalities \eqref{3d-h1}, \eqref{oh-v2}, and 
\eqref{oh-e1} imply that the 
sequence of traces 
$\tr^\fc(S_3(K)) = S_2(\fc)$ for any face $\fc\in\FC(K)$ is 
a compatible exact sequence. Hence, $S_3(K)$ is a compatible
exact sequence.

Finally, invoking \cite[Corollary 3.2]{ChristiansenGillette15},
we get the minimality of the sequence $S_3(K)$ by the simple observation that  
\[
\oo{H_3}= \oo{H_3^g}, \quad
\oo{E_3}= \oo{E_3^g},\quad
\oo{V_3}= \oo{V_3^g},\quad
\oo{W_3}= \oo{W_3^g}.
\]

This completes the proof of Theorem 
\ref{thm:3d-general}. 
\end{proof}

\section{Proofs of the applications in Section \ref{sec:applications}}
\label{sec:proof2}
In this section, we prove all the results of 
Section \ref{sec:applications} by applying Theorem \ref{thm:1d-general} (for one dimension),
Theorem \ref{thm:2d-general} (for two dimensions),
and Theorem \ref{thm:3d-general} (for three dimensions).

\subsection{The one-dimensionanl case}
The proof of Theorem \ref{thm:1d} is a direct application of 
Theorem \ref{thm:1d-general}.

\subsection{The two-dimensional case}

We first present the following result 
on  exact sequences on the whole space $\bR^2$. \gfu{We give its proof in Appendix B.}
\begin{lemma}
 \label{lemma:2d-exact}
The following four sequences on $\bR^2$ with 
coordinates $(x,y)$
 are exact for $k\ge 0$:
 
 {\centering
\begin{tabular}{l c c c c c}
$\mathrm{S_{1,k}^{
2d}}:$ &
$ \pol_{k+2}$&
$\overset{\grads}{\longrightarrow}$ &
$\bld\pol_{k+1} $&
$\overset{\curls}{\longrightarrow}$ &
$\pol_{k}$,\\
\\
 $\mathrm{S_{2,k}^{2d}
 }$: &$\pol_{k+1}$&
$\overset{\grads}{\longrightarrow}$ &
$\bld\pol_{k}\oplus\bld x
\times\widetilde\pol_{k}$
&
$\overset{\curls}{\longrightarrow}$ &
$\pol_{k}$,\\
\\
$\mathrm{S_{3,k}^{
2d}}$: &$\qol_{k}\oplus \{x^{k+1},y^{k+1}\}$&
$\overset{\grads}{\longrightarrow}$ &
$\bld\qol_{k}\oplus\bld x
\times\{x^ky^k\}$
&
$\overset{\curls}{\longrightarrow}$ &
$\qol_{k}$,\\
 \\
$\mathrm{S_{4,k}^{
2d
}}$: &$\qol_{k+1}$&
$\overset{\grads}{\longrightarrow}$ &
$\bld\qol_{k}\oplus
\left(
\begin{tabular}{c}
 $y^{k+1}\pol_k(x)$\\
 $x^{k+1}\pol_k(y)$
\end{tabular}
\right)$
&
$\overset{\curls}{\longrightarrow}$ &
$\qol_{k}$.
\end{tabular}
}\\
\end{lemma}

The third sequence $\mathrm{S_{2,k}^{3d}}$
is new to the best of the 
authors' knowledge; the other three are well-known. 
Note that all these four sequences has good  ``symmetry'' 
in the sense that they are invariant under the 
coordinate permutation $(x,y)\longrightarrow (y,x)$. 

Now, we are ready to prove the two-dimensional 
results of Theorem \ref{thm:2d-t}, Theorem \ref{thm:2d-s},
and Theorem \ref{thm:2d-p} by
applying Theorem \ref{thm:2d-general}.

\subsubsection*{Proof of Theorem \ref{thm:2d-t}}
\begin{proof}
Let us fit the sequences $\satri$ and 
$\sbtri$ into the framework of 
Theorem \ref{thm:2d-general}.

For the first sequence $\satri$, we have that the 
set of complete trace sequences is
\[
S_1(\dK):=
\{
S_1(e):\;\;\;
\pol_{k+2}(e)\longrightarrow 
\pol_{k+1}(e)\;\forall e\in \EG(K),
\}
\]
and the $S_1(\dK)$-admissible sequence is 
$\satri$ itself. We also have 
\[
\dim \delta H_2^g = 3(k+2)+k(k+1)/2-(k+3)(k+4)/2 = 0,
\]
hence $\delta H_2^g(K)=\{0\}$. 
So, $\satri$ is a commuting exact sequence.
Similarly, we conclude that $\sbtri$
is also a commuting exact sequence.
This completes the proof of Theorem \ref{thm:2d-t}.
\end{proof}

\subsubsection*{Proof of Theorem \ref{thm:2d-s}}
\begin{proof}
Let us fit the sequences $\sasqr$, $\sbsqr$, 
$\scsqr$, and 
$\sdsqr$ into the framework of 
Theorem \ref{thm:2d-general}.

For the first sequence $\sasqr$, we have that the 
set of complete trace sequences is
\[
S_1(\dK):=
\{
S_1(e):\;\;\;
\pol_{k+2}(e)\longrightarrow 
\pol_{k+1}(e)\;\forall e\in \EG(K),
\}
\]
and the $S_1(\dK)$-admissible sequence is 
$\mathrm{S_{1,k}^{2d}}$ in Lemma \ref{lemma:2d-exact}.
We also have
the space
$
\delta H_2^g:= \delta H_{k+2}^{2,I}
$
satisfy all the four properties 
of $\delta H_2^g$ in Theorem \ref{thm:2d-general}, 
in particular, we have 
\[
 \dim \delta H_2^g = 4(k+2) + (k-1)k/2-(k+3)(k+4)/2=2.
\]
Hence, $\sasqr$ is a commuting exact sequence.

The proof for the second sequence $\sbsqr$
is identical to that for the first one.

For the third sequence $\scsqr$, we have that the 
set of complete trace sequences is
\[
S_1(\dK):=
\{
S_1(e):\;\;\;
\pol_{k+1}(e)\longrightarrow 
\pol_{k}(e)\;\forall e\in \EG(K),
\}
\]
and the $S_1(\dK)$-admissible sequence is 
$\mathrm{S_{3,k}^{2d}}$ in Lemma \ref{lemma:2d-exact}.
We also have
the space
$
\delta H_2^g:= \delta H_{k+1}^{2,I}
$
satisfy all the four properties 
of $\delta H_2^g$ in Theorem \ref{thm:2d-general}, 
in particular, we have 
\[
 \dim \delta H_2^g = 4(k+1) + (k-1)^2-\left(
 (k+1)^2+2\right)=2.
\]
Hence, $\scsqr$ is a commuting exact sequence.

For the last sequence $\sdsqr$, we have the 
set of complete trace sequences is
\[
S_1(\dK):=
\{
S_1(e):\;\;\;
\pol_{k+1}(e)\longrightarrow 
\pol_{k}(e)\;\forall e\in \EG(K),
\}
\]
and the $S_1(\dK)$-admissible sequence is 
$\mathrm{S_{4,k}^{2d}}$ in Lemma \ref{lemma:2d-exact}.
We also have
the space
$
\delta H_2^g:= \{0\}
$
satisfy all the four properties 
of $\delta H_2^g$ in Theorem \ref{thm:2d-general}, 
since 
\[
 \dim \delta H_2^g = 4(k+1) + k^2-
 (k+2)^2=0.
\]
Hence, $\sdsqr$ is also a commuting exact sequence.
This completes the proof of Theorem \ref{thm:2d-s}.
\end{proof}

\subsubsection*{Proof of Theorem \ref{thm:2d-p}}
\begin{proof}
Let us fit the sequences 
$\mathrm{S_{1,k}^{poly}}$ and
$\mathrm{S_{2,k}^{poly}}$
into the framework of 
Theorem \ref{thm:2d-general}.

For the first sequence 
$\mathrm{S_{1,k}^{poly}}$, we have the 
set of complete trace sequences is
\[
S_1(\dK):=
\{
S_1(e):\;\;\;
\pol_{k+2}(e)\longrightarrow 
\pol_{k+1}(e)\;\forall e\in \EG(K),
\}
\]
and the $S_1(\dK)$-admissible sequence is 
$\mathrm{S_{1,k}^{2d}}$ in Lemma \ref{lemma:2d-exact}.
The proof of the space $\delta H_2^g:=\delta H_{k+2}^{2,II}$
satisfying all the four properties of $\delta H_2^g$ 
in Theorem \ref{thm:2d-general}
is not trival. It is given in the proof of 
\cite[Theorem 2.6]{CockburnFuM2D}.

The proof for the second sequence
is identical to that for the first one.
This completes the proof of Theorem \ref{thm:2d-p}.
\end{proof}

\subsection{The three-dimensional case}
As we did for the two-dimensional case, we first present 
the following result 
on  exact sequences on the whole space $\bR^3$.
\gfu{We give its proof in Appendix B.}
\begin{lemma}
 \label{lemma:3d-exact}
The following six sequences on $\bR^3$
 are exact for $k\ge 0$.
The first two sequences
 are the famous polynomial de Rham sequences that
 contain
 polynomials of certain degree:

\

 {
\begin{tabular}{l c c c c c c c}
$\mathrm{S_{1,k}^{
3d}:}$ &
$ \pol_{k+3}$&
$\overset{\grads}{\longrightarrow}$ &
$\bld\pol_{k+2} $&
$\overset{\curls}{\longrightarrow}$ &
$\bpol_{k+1}$
&
$\overset{\divs}{\longrightarrow}$ &
$\pol_{k}$,\\
 $\mathrm{S_{2,k}^{3d}
 }:$ &$\pol_{k+1}$&
$\overset{\grads}{\longrightarrow}$ &
$\bld\pol_{k}\oplus\bld x
\times\bld{\widetilde\pol}_{k}$
&
$\overset{\curls}{\longrightarrow}$ &
$\bld\pol_{k}\oplus\bld x
{\widetilde\pol}_{k}$
&
$\overset{\divs}{\longrightarrow}$ &
$\pol_{k}$,\\
\end{tabular}

The next two sequences have spaces containing tensor product polynomials of certain degree.

\

\resizebox{0.95\columnwidth}{!}
{
\begin{tabular}{l c c c c c c c}
$\mathrm{S_{3,k}^{
3d}:}$ &
$ \qol_{k}
$
&$\overset{\grads}{\longrightarrow}$ &
$ \bqol_{k}
$&
$\overset{\curls}{\longrightarrow}$ &
$\bld\qol_{k}
$&
$\overset{\divs}{\longrightarrow}$ &
$\qol_{k}
$,\\
&
$ \oplus \left\{
\begin{tabular}{c}
 $x^{k+1},$\\$y^{k+1},$
 \\
 $z^{k+1}$
\end{tabular}
\right\}
$& &
$ \oplus \bld x\times 
\left\{
\begin{tabular}{c}
$ y^kz^k\grads x$, \\
$z^kx^k\grads y,$\\
$x^ky^k\grads z$
\end{tabular}
\right\}
$&
&
$\oplus \bld x\,\{x^ky^kz^k\}
$&
 &\\
$\mathrm{S_{4,k}^{3d}:}$ &
$ \qol_{k+1}$&
$\overset{\grads}{\longrightarrow}$ &
$ 
\left(
\begin{tabular}{c}
$ \pol_{k,k+1,k+1}$ \\
$\pol_{k+1,k,k+1}$\\
$\pol_{k+1,k+1,k}$
\end{tabular}
\right)
$&
$\overset{\curls}{\longrightarrow}$ &
$\left(
\begin{tabular}{c}
$ \pol_{k+1,k,k}$ \\
$\pol_{k,k+1,k}$\\
$\pol_{k,k,k+1}$
\end{tabular}
\right)
$&
$\overset{\divs}{\longrightarrow}$ &
$\qol_{k}
$.\\
\end{tabular}
}

The last two sequences 
contain polynomials of certain degree in the $(x,y)-$plane, and 
have some spaces with tensor product structure in the $(x,z)-$ and $(y,z)-$plane.

\

\resizebox{0.95\columnwidth}{!}
{
\begin{tabular}{l c c c c c c c}
 $\mathrm{S_{5,k}^{3d}
}:$ &
 $ \pol_{k|k}
 $&
 $\overset{\grads}{\longrightarrow}$ &
 $ \bpol_{k|k}
$
&
$\overset{\curls}{\longrightarrow}$ &
$\bld\pol_{k|k}$&
$\overset{\divs}{\longrightarrow}$ &
$\pol_{k|k}
$,\\
&\!\!\!\!\!\!\!\!
 $\oplus\widetilde\pol_{k+1}(x,y)\oplus
\{ z^{k+1}\}
 $
 &
 &
 $\oplus 
 \left(\!\!\!\!
 \begin{tabular}{c}
$  \left(\!\!\!\!
  \begin{tabular}{c}
 $y$\\
 $-x$
 \end{tabular}
\!\!\!\! \right)\widetilde\pol_{k}(x,y)$\\
$\widetilde\pol_{k+1}(x,y)\,z^k$
\end{tabular}
\!\!\!\! \right)
$
&
&
$\oplus 
 \left(\!\!\!\!
 \begin{tabular}{c}
$0$\\
$0$\\
$\widetilde\pol_{k}(x,y)\,z^{k+1}$
\end{tabular}
\!\!\!\! \right)
$&
\vspace{.1cm}
\\
$\mathrm{S_{6,k}^{3d}
}:$ &
$ \pol_{k+1|k+1}$&
$\overset{\grads}{\longrightarrow}$ 
&
$ 
 \left(\!\!\!\!
 \begin{tabular}{c}
$
\bld{N}_{k}(x,y)\otimes\pol_{k+1}(z)$\\	
$\pol_{k+1|k}$
\end{tabular}
\!\!\!\! \right)
$&
$\overset{\curls}{\longrightarrow}$ &
$\left(\!\!\!\!
 \begin{tabular}{c}
$ \bld{{RT}}_{k}(x,y)\otimes\pol_{k}(z)$\\	
$\pol_{k|k+1}$
\end{tabular}
\!\!\!\! \right)
$&
$\overset{\divs}{\longrightarrow}$ &
$\pol_{k|k}
$.\\
\end{tabular}
}}\\
\end{lemma}

The third sequence $\mathrm{S_{3,k}^{3d}}$ and 
the fifth sequence $\mathrm{S_{5,k}^{3d}}$ are new to the best of the authors'
knowledge; the other four are well-known. 
Note that all these six sequences, except the fifth, 
have good ``symmetry'' in the sense that they are 
invariant under any coordinate permutation. The fifth sequence is only 
invariant under the coordinate permutation $(x,y,z)\longrightarrow (y,x,z)$.

To further simplify the notation, 
we use the so-called $M$-index introduced in \cite{CockburnFuSayas16} and used in 
\cite{CockburnFuM3D} to obtain various finite element spaces 
admitting $M$-decompositions in three dimensions.
The definition of an $M$-index is given as follows.
\begin{definition}[The $M$-index]
The $M$-index of the space $V(K)\times W(K)\subset H(\mathrm{div},K)\times H^1(K)$ is 
the number 
\begin{alignat*}{2}
I_{M}(V\times W):=&\;\dim M(\dK) &&-\dim\tr_V\{v\in V(K):\, \divs v=0\}
                                     \\
                                     &&&-\dim \tr_H\{w\in W(K):\, \grads w=0\},
\end{alignat*}
where 
\[
M(\dK)=
\{
\mu\in L^2(\dK):\;
\mu|_{\fc}\in M(\fc)\quad \forall \fc\in\FC(K)
\}
\]
is a finite element space defined on the boundary
$\dK$ of a polyhedron $K$.
\end{definition}

Using the definition of an $M$-index, we have 
the number in the 
right hand side of property $\mathrm{(iv)}$ of $\delta E_3^g(K)$ 
in Theorem \ref{thm:3d-general}
is 
nothing but 
$I_M(V_3^g\times W_3^g)$
with the trace space 
\begin{align}
\label{trace-space}
  M(\dK):=
 \{
\mu\in L^2(\dK):\;
\mu|_{\fc}\in W_2(\fc)\quad \forall \fc\in\FC(K) 
 \}.
\end{align}
That is, 
\begin{align}
\label{m-index}
I_M(V_3^g\times W_3^g) =&\; \sum_{\fc\in \FC(K)}W_2(\fc)
 +\dim \oo{W_3^g}(K)\\
 &\;
 +\dim \left\{v\in \oo{V_3^g}(K):\divs v=0\right\}-\dim V_3^g(K).
 \nonumber
\end{align}

To see this, we have 
\begin{align*}
& \sum_{\fc\in \FC(K)}W_2(\fc)
 +\dim \oo{W_3^g}(K)
 +\dim \left\{v\in \oo{V_3^g}(K):\divs v=0\right\}-\dim V_3^g(K)\;\;\\
 &=\dim M(\dK) +\dim W_3^g(K)-1 +\dim \left\{
 v\in V_3^g(K):\;\divs v=0
 \right\}\\
 &\hspace{1cm}-
 \dim \tr_V\left\{
 v\in V_3^g(K):\;\divs v=0
 \right\} - \dim V_3^g(K)\\
& =\dim M(\dK) +\dim W_3^g(K)-1 - \dim \tr_V\left\{
 v\in V_3^g(K):\;\divs v=0
 \right\}\\
 &\hspace{1cm}-
 \dim \divs V_3^g(K)\\
&  =\dim M(\dK) -
\dim \tr_H\left\{
 w\in W_3^g(K):\;\grads w=0
 \right\}\\
 &\hspace{2.28cm}- \dim \tr_V\left\{
 v\in V_3^g(K):\;\divs v=0
 \right\}\\
 &=I_M(V_3^g\times W_3^g).
\end{align*}
From now on, the trace space $M(\dK)$ will always be of the form \eqref{trace-space} where 
the spaces $W_2(\fc)$ on each face vary in different locations.

Since the computation of
the $M$-index for various polynomial spaces on the four reference polyhedra was given in
\cite{CockburnFuM3D}, we can directly use those results to verify property 
$\mathrm{(iv)}$
of $\delta E_3^g(K)$.

Now, we are ready to prove the three-dimensional results of 
Theorem \ref{thm:3d-t} to Theorem \ref{thm:3d-pyramid} by applying the 
general result of Theorem \ref{thm:3d-general}.
\subsubsection*{Proof of Theorem \ref{thm:3d-t}}
\begin{proof}
Let us fit the sequences $\satet$ and $\sbtet$ into the framework of Theorem 
\ref{thm:3d-general}. Here the element 
\[
 K =
 \{(x,y,z): 0<x <1,\; 0<y<1,\; 0<z<1,\; x+y+z<1\}
\]
is the reference tetrahedron with four square faces, 
six edges, and four vertices.

For the first sequence $\satet$, we have the set of complete trace sequences
is 
\[
 S_2(\dK):=\{\mathrm{\satrip}(\fc):\;\;\forall \fc\in\FC(K)\},
\]
and the given $S_2(\dK)$-admissible sequence is $\mathrm{S_{1,k}^{3d}}$
in Lemma \ref{lemma:3d-exact}.
We also have 
\begin{align*}
  \dim \delta H_3^g = &\; 4 + 6(k+2) + 4\,(k+1)(k+2)/2+k(k+1)(k+2)/6\\
  &\;-
(k+4)(k+5)(k+6)/6\\
=&\; 0 \\
  \dim \delta E_3^g = &\; I_M(V_3^g\times W_3^g)\\
=&\; 0,
\end{align*}
hence $\delta H_3^g =\{0\}$, and 
$\delta E_3^g = \{0\}$. So, $\satet$ is a commuting exact sequence.
Similarly, we conclude that $\sbtet$ is also a commuting exact sequence.
This completes the proof of Theorem \ref{thm:3d-t}.
\end{proof}

\subsubsection*{Proof of Theorem \ref{thm:3d-hex}}
\begin{proof}
Let us fit the sequences $\sacub$, $\sbcub$, $\sccub$, and $\sdcub$ 
into the framework of Theorem 
\ref{thm:3d-general}.
Here the element 
\[
 K =
 \{(x,y,z): 0<x <1,\; 0<y<1,\; 0<z<1\}
\]
is the reference cube with six square faces, 
twelve edges, and eight vertices.

For the first sequence $\sacub$, we have the set of complete trace sequences
is 
\[
 S_2(\dK):=\{\mathrm{\sasqrp}(\fc):\;\;\forall \fc\in\FC(K)\},
\]
and the given $S_2(\dK)$-admissible sequence is $\mathrm{S_{1,k}^{3d}}$
in Lemma \ref{lemma:3d-exact}.
It is then easy to verify that 
$\delta H_3^g:=\delta H_{k+3}^{3,I}$ satisfies all the four properties of 
$\delta H_3^g$, and
$\delta E_3^g:=\delta E_{k+2}^{3,I}$ satisfies all the four properties of 
$\delta E_3^g$ in Theorem \ref{thm:3d-general}.
In particular, we have 
\begin{align*}
  \dim \delta H_3^g = &\; 8 + 12(k+2) + 6\,k(k+1)/2+(k-2)(k-1)k/6\\
  &\;-
(k+4)(k+5)(k+6)/6\\
=&\; 3(k+4), \\
  \dim \delta E_3^g = &\; I_M(V_3^g\times W_3^g)
=\; 3(k+2),
\end{align*}
where the last equality is due to \cite[Theorem 2.12]{CockburnFuM3D}.

For the sake of completeness, here we 
we present a proof of property (iii) for
$\delta E_3^g(K)$, 
which is a bit more difficult to verify that the other properties.
Given a function $p\in \widetilde\pol_{k+1}(y,z)$, we have 
\[
\curls\left( x\,p\,(y\grads z-z\grads y)\right) = (k+3)p\,x\grads x-p\,y\grads y-p\,z\grads z.
\]
This means that the function $\curls\left( x\,p\,(y\grads z-z\grads y)\right)
\in \curls \delta E_3^g(K)$ has normal trace equal
to {\it zero} on the three faces $x=0$, $y=0$, and $z=0$ of the unit cube $K$, and has normal trace equals to 
$(k+3)p\in \widetilde\pol_{k+1}(y,z)$ on the face $x=1$. 
Similar results hold for 
\begin{align*}
\curls\left( y\,q\,(z\grads x-x\grads z)\right) = -q\,x\grads x+(k+3)q\,y\grads y-q\,z\grads z,\\
\curls\left( z\,r\,(x\grads y-y\grads x)\right) = -r\,x\grads x-r\,y\grads y+(k+3)r\,z\grads z,
\end{align*}
with $q\in \widetilde\pol_{k+1}(z,x)$, and $r\in \widetilde\pol_{k+1}(x,y)$.
Using this fact, we have any function in $\curls \delta E_3^g(K)$ 
has a normal trace equal to {\it zero} on one face and be a polynomial of degree 
$k+1$ on its opposite (parallel) face for at least one pair of parallel faces. 
On the other hand, if a function in $V_3^g(K)=\bpol_{k+1}$ has normal trace 
equal to {\it zero} on one face, the normal trace on its opposite face must be a polynomial of degree no greater
than $k$. 
Hence, $\tr_V\curls \delta E_3^g(K)\cap \tr_V V_3^g(K)=\{0\}$. This implies 
property $\mathrm{(iii)}$ of $\delta E_3^g(K)$.

So, $\sacub$ is a commuting exact sequence.

The proof for the second sequence $\sbcub$ is identical to that for the 
first one.

For the third sequence $\sccub$, we have the set of complete trace sequences
is 
\[
 S_2(\dK):=\{\mathrm{\scsqr}(\fc):\;\;\forall \fc\in\FC(K)\},
\]
and the given $S_2(\dK)$-admissible sequence is $\mathrm{S_{3,k}^{3d}}$
in Lemma \ref{lemma:3d-exact}.
It is then trivial to verify that
$\delta H_3^g:=\delta H_{k+1}^{3,II}$ satisfies all the four properties of 
$\delta H_3^g$, and
$\delta E_3^g:=\delta E_{k+1}^{3,II}$ satisfies all the four properties of 
$\delta E_3^g$ in Theorem \ref{thm:3d-general}.
In particular, we have, for $k\ge 1$,
\begin{align*}
  \dim \delta H_3^g = &\; 
   8 + 12 k + 6(k-1)^2 + {(k-1)^3}-\left({(k+1)^3+3}\right) \\
   =&\; 9,\\
  \dim \delta E_3^g = &\; I_M(V_3^g\times W_3^g)
=\; 6,
\end{align*}
where the last equality is due to \cite[Theorem 2.11]{CockburnFuM3D}.
So, $\sccub$ is a commuting exact sequence.

For the last sequence $\sdcub$, we have the set of complete trace sequences
is 
\[
 S_2(\dK):=\{\mathrm{\sdsqr}(\fc):\;\;\forall \fc\in\FC(K)\},
\]
and the given $S_2(\dK)$-admissible sequence is $\mathrm{S_{4,k}^{3d}}$
in Lemma \ref{lemma:3d-exact}.
We also have the space 
$\delta H_3^g:=\{0\}$ satisfies all the four properties of 
$\delta H_3^g$, and
$\delta E_3^g:=\{0\}$ satisfies all the four properties of 
$\delta E_3^g$ in Theorem \ref{thm:3d-general} since
\begin{align*}
  \dim \delta H_3^g = &\; 
   8 + 12 k + 6 k^2 + {k^3}-{(k+2)^3} \\
   =&\; 0,\\
  \dim \delta E_3^g = &\; I_M(V_3^g\times W_3^g)
=\; 0.
\end{align*}
So, $\sdcub$ is a commuting exact sequence.
This completes the proof of Theorem \ref{thm:3d-hex}.
\end{proof}

\subsubsection*{Proof of Theorem \ref{thm:3d-prism}}
\begin{proof}
Let us fit the sequences $\sapri$, $\sbpri$, $\scpri$, and $\sdpri$ 
into the framework of Theorem 
\ref{thm:3d-general}.
Here the element 
\[
 K =
 \{(x,y,z): 0<x ,\; 0<y,\; 0<z<1,\; x+y<1 \}
\]
is the reference prism with five faces (two triangular faces
and three square faces), 
nine edges, and six vertices.

For the first sequence $\sapri$, we have the set of complete trace sequences
is 
\[
 S_2(\dK):=\left\{S_2(\fc):\;\;
 \begin{tabular}{l l}
$S_2(\fc)=\satrip(\fc)$ & if $\fc$ is a triangle,\\
$S_2(\fc)=\sasqrp(\fc)$ & if $\fc$ is a square
 \end{tabular}
\right\},
\]
and the given $S_2(\dK)$-admissible sequence is $\mathrm{S_{1,k}^{3d}}$
in Lemma \ref{lemma:3d-exact}.
It is then trivial to verify that 
$\delta H_3^g:=\delta H_{k+3}^{3,III}$ satisfies all the four properties of 
$\delta H_3^g$, and
$\delta E_3^g:=\delta E_{k+2}^{3,III}$ satisfies all the four properties of 
$\delta E_3^g$ in Theorem \ref{thm:3d-general}.
In particular, we have 
\begin{align*}
  \dim \delta H_3^g = &\; 
  6 + 9 (k+2) +(k+1)(k+2) + 3\,k(k+1)/2 \\
&\;+ {(k-1)k(k+1)}/{6}-
 {(k+4)(k+5)(k+6)}/{6}\\
 =&\; k+6,\\
  \dim \delta E_3^g = &\; I_M(V_3^g\times W_3^g)
=\; k+3,
\end{align*}
where the last equality is due to \cite[Theorem 2.8]{CockburnFuM3D}.
So, $\sapri$ is a commuting exact sequence.

The proof for the second sequence $\sbpri$ is identical to that for the 
first one.

For the third sequence $\scpri$, we have the set of complete trace sequences
is 
\[
 S_2(\dK):=\left\{S_2(\fc):\;\;
 \begin{tabular}{l l}
$S_2(\fc)=\sbtri(\fc)$ & if $\fc$ is a triangle,\\
$S_2(\fc)=\widetilde{\scsqr}(\fc)$ & if $\fc$ is a square
 \end{tabular}
\right\},
\]
and the given $S_2(\dK)$-admissible sequence is $\mathrm{S_{5,k}^{3d}}$
in Lemma \ref{lemma:3d-exact}.
It is then trivial to verify that
$\delta H_3^g:=\delta H_{k+1}^{3,III}$ satisfies all the four properties of 
$\delta H_3^g$, and
$\delta E_3^g:=\delta E_{k+1}^{3,III}$ satisfies all the four properties of 
$\delta E_3^g$ in Theorem \ref{thm:3d-general}.
In particular, we have, for $k\ge 1$,
\begin{align*}
  \dim \delta H_3^g = &\; 
   6 + 9 k + (k-1)k + {3(k-1)^2}+(k-2)(k-1)^2/2\\
   &\;  -\left({(k+1)^2(k+2)/2+k+3}\right) \\
   =&\; k+4,\\
  \dim \delta E_3^g = &\; I_M(V_3^g\times W_3^g)
=\; k+2,
\end{align*}
where the last equality is due to \cite[Theorem 2.7]{CockburnFuM3D}.
So, $\scpri$ is a commuting exact sequence.

For the last sequence $\sdpri$, we have the set of complete trace sequences
is 
\[
 S_2(\dK):=\left\{S_2(\fc):\;\;
 \begin{tabular}{l l}
$S_2(\fc)=\sbtri(\fc)$ & if $\fc$ is a triangle,\\
$S_2(\fc)={\sdsqr}(\fc)$ & if $\fc$ is a square
 \end{tabular}
\right\},
\]
and the given $S_2(\dK)$-admissible sequence is $\mathrm{S_{6,k}^{3d}}$
in Lemma \ref{lemma:3d-exact}.
We also have the space 
$\delta H_3^g:=\{0\}$ satisfies all the four properties of 
$\delta H_3^g$, and
$\delta E_3^g:=\{0\}$ satisfies all the four properties of 
$\delta E_3^g$ in Theorem \ref{thm:3d-general} since
\begin{align*}
  \dim \delta H_3^g = &\; 
   6 + 9 k + (k-1)k + {3\,k^2}+(k-1)k^2/2\\
   &\;  -{(k+2)^2(k+3)/2} \\
   =&\; 0,\\
  \dim \delta E_3^g = &\; I_M(V_3^g\times W_3^g)
=\; 0.
\end{align*}
So, $\sdpri$ is a commuting exact sequence.
This completes the proof of Theorem \ref{thm:3d-prism}.
\end{proof}

\subsubsection*{Proof of Theorem \ref{thm:3d-pyramid}}
\begin{proof}
Let us fit the sequences $\sapyr$, $\sbpyr$, $\scpyr$, and $\sdpyr$ 
into the framework of Theorem 
\ref{thm:3d-general}.
Here the element 
\[
 K =
 \{(x,y,z): 0<x<1 ,\; 0<y<1,\; 0<z<1,\; x+z<1,\;
 y+z<1\}
\]
is the reference prism with five faces (four triangular faces
and one square face), 
eight edges, and five vertices.

For the first sequence $\sapyr$, we have the set of complete trace sequences
is 
\[
 S_2(\dK):=\left\{S_2(\fc):\;\;
 \begin{tabular}{l l}
$S_2(\fc)=\satrip(\fc)$ & if $\fc$ is a triangle,\\
$S_2(\fc)=\sasqrp(\fc)$ & if $\fc$ is a square
 \end{tabular}
\right\},
\]
and the given $S_2(\dK)$-admissible sequence is $\mathrm{S_{1,k}^{3d}}$
in Lemma \ref{lemma:3d-exact}.
It is then trivial to verify that 
$\delta H_3^g:=\delta H_{k+3}^{3,IV}$ satisfies all the four properties of 
$\delta H_3^g$, and
$\delta E_3^g:=\delta E_{k+2}^{3,IV}$ satisfies all the four properties of 
$\delta E_3^g$ in Theorem \ref{thm:3d-general}.
In particular, we have 
\begin{align*}
  \dim \delta H_3^g = &\; 
  5 + 8 (k+2) +2(k+1)(k+2) + k(k+1)/2 \\
&\;+ {(k-1)k(k+1)}/{6}-
 {(k+4)(k+5)(k+6)}/{6}\\
 =&\; 2\,k+5,\\
  \dim \delta E_3^g = &\; I_M(V_3^g\times W_3^g)
=\; 3,
\end{align*}
where the last equality is due to \cite[Theorem 2.6]{CockburnFuM3D}.
So, $\sapyr$ is a commuting exact sequence.

The proofs for the other three  sequences are similar to that for the 
first one, and hence are omitted.
This completes the proof of Theorem \ref{thm:3d-pyramid}.
\end{proof}

\section{Conclusion}
\label{sec:conclude}

We presented a systematic construction of  commuting exact sequences
on polygonal/polyhedral elements. 
The systematic construction is applied to the reference triangle, reference square, and 
general polygon
in two dimensions, 
and the reference tetrahedron, reference cube, reference prism, 
and reference pyramid in three dimensions to obtain concrete 
commuting exact sequences. 
We obtain {\it mapped} commuting exact sequences on a two-dimensional hybrid mesh
consists of triangles and parallelograms, and 
{\it non-mapped} 
commuting exact sequences on a general 
two-dimensional polygonal mesh.
We also obtain {\it mapped} commuting exact sequences on a 
three-dimensional hybrid mesh
with elements obtained via affine mapping from one of
the four reference polyhedra. The actual implementation of shape functions for 
our four family of commuting exact sequences on hybrid polyhedral meshes constitutes the subject of ongoing work. 

As already pointed out in \cite{CockburnFuM3D}, 
finding stable pairs of finite element spaces  
defining mixed methods for the diffusion problem
on general polyhedral elements is a very complicated task, mainly due to 
the need of characterizing the  {\it solenoidal bubble space } 
\[
 \{v\in \oo{V_3^g}(K):\;\;\divs v=0\},
\]
that is, the subspace of $V_3^g$ containing divergence-free
functions with zero normal trace on the boundary.
Finding exact sequences with a commuting diagram property shares exactly the same difficulty.
However, Theorem \ref{thm:3d-general} 
does shed light on a promising approach to carry out the construction.
In particular, we can
take the given  exact sequence to be 
$\mathrm{S_{1,k}^{3d}}$, which is 
a $S_2(\dK)$-admissible exact sequence for 
the set of commuting trace sequences
\[
 S_2(\dK):=\{\mathrm{S_{1,k}^{poly}}(\fc):\;\;\forall \fc\in\FC(K)
 \},
\]
or to be
$\mathrm{S_{2,k}^{3d}}$, which is 
a $S_2(\dK)$-admissible exact sequence for 
the set of commuting trace sequences
\[
 S_2(\dK):=\{\mathrm{S_{2,k}^{poly}}(\fc):\;\;\forall \fc\in\FC(K)
 \}.
\]
Then, the only task left
is to find the spaces $\delta H_3^g(K)$ and $\delta E_3^g(K)$
that satisfy the properties described in Theorem \ref{thm:3d-general}.
The actual construction of commuting exact sequences on a 
general polyhedron is currently under way.

\section*{Acknowledgements}
The authors would like to thank Douglas Arnold from the University of Minnesota for many helpful discussions on commuting exact sequences.

\section*{Appendix A: the harmonic interpolators}
Here we use our notation to rewrite the harmonic interpolators (in one-, two- and  three-space dimensions) introduced in \cite{ChristiansenMuntheKaasOwren11}. To simplify the notation, we let the domain of these interpolations be space of smooth fields. 
However, we only need these spaces to be regular enough such that the 
related differential and traces operators used in the harmonic interpolators make sense; see \cite[Section 2.3]{ChristiansenRapetti16} for the related spaces 
with minimal regularity. 

\subsection*{One-dimensional harmonic interpolators}
Let $K\in \bR$ be a segment. We denote $\Pi_H^1$ and $\Pi_W^1$ 
being the harmonic interpolators related to $H^1$- and $L^2$-fields.
We denote 
\[
 H_1(K)\longrightarrow W_1(K)
\]
as the related compatible exact sequence.

The harmonic interpolators 
$ \Pi_H^1\times \Pi_W^1:\;\HH(\bar{K})\times \HH(\bar{K})\rightarrow\; H_1(K)\times W_1(K)$ 
are the following ones 
satisfying
\begin{alignat*}{2}
 \Pi_H^1 u(\vt) =&\; u(\vt)&&\quad \forall \vt\in \VT(K),\\
 (\grads\Pi_H^1 u,v)_K =&\; (\grads u,v)_K&&\quad 
 \forall v\in \grads\overset{\circ}{H_1}(K)
\end{alignat*}
and 
\begin{alignat*}{2}
 (\Pi_W^1 u,v)_K =&\; (u,v)_K&&\quad \forall v\in \oo{W_1}(K)\oplus \pol_0(K).
\end{alignat*}

\subsection*{Two-dimensional harmonic interpolators}
Let $K\in \bR^2$ be a segment. We denote $\Pi_H^2$, $\Pi_E^2$ and $\Pi_W^2$ 
being the harmonic interpolators related to $H^1$-, $H(\mathrm{curl})$-, and $L^2$-fields.
We denote 
\[
 H_2(K)\longrightarrow E_2(K)\longrightarrow W_2(K)
\]
as the related compatible exact sequence.

The harmonic interpolators 
$ \Pi_H^2\times \Pi_E^2\times \Pi_W^2:\;\HH(\bar{K})\times \bld\HH(\bar{K})\times \HH(\bar{K})\rightarrow\; 
H_2(K)\times E_2(K)\times W_2(K)$ 
are the following ones 
satisfying
\begin{alignat*}{2}
 \Pi_H^2 u(\vt) =&\; u(\vt)&&\quad \forall \vt\in \VT(K),\\
 \left(\grads \tr_H^\eg(\Pi_H^2 u),v
 \right)_\eg =&\; (\grads \tr_H^\eg(u),v)_\eg&&\quad \forall v\in 
 \grads\overset{\circ}{H_1}(\eg),
 \forall \eg\in \EG(K),\\
 (\grads\Pi_H^2 u,v)_K =&\; (\grads u,v)_K&&\quad \forall 
 v\in\grads \overset{\circ}{H_2}(K),
\end{alignat*}
and
\begin{alignat*}{2}
\left(\tr_E^\eg(\Pi_E^2 u),v\right)_\eg =&\; \left(\tr_E^\eg(u),v\right)_\eg&&
\quad \forall v\in \oo{W_1}(\eg)\oplus \pol_0(\eg),
 \forall \eg\in \EG(K),\\
 (\curls \Pi_E^2 u,v)_K =&\; (\curls u,v)_K&&
 \quad \forall v\in \curls\overset{\circ}{E_2}(K),\\
  (\Pi_E^2 u,v)_K =&\; (u,v)_K&&
 \quad \forall v\in \{v\in \overset{\circ}{E_2}(K):\;\curls v=0\},
\end{alignat*}
and 
\begin{alignat*}{2}
 (\Pi_W^2 u,v)_K =&\; (u,v)_K&&\quad \forall v\in \oo{W_2}(K)\oplus \pol_0(K).
\end{alignat*}

\subsection*{Three-dimensional harmonic interpolators}
Let $K\in \bR^3$ be a polyhedron. We denote $\Pi_H^3$, $\Pi_E^3$, $\Pi_V^3$ and $\Pi_W^3$ 
being the harmonic interpolators related to $H^1$-, $H(\mathrm{curl})$-, $H(\mathrm{div})$- and $L^2$-fields.
We denote 
\[
 H_3(K)\longrightarrow E_3(K)\longrightarrow V_3(K)\longrightarrow W_3(K)
\]
as the related compatible exact sequence.

The harmonic interpolators 
$ \Pi_H^3\times \Pi_E^3\times \Pi_V^3\times \Pi_W^3:\;\HH(\bar{K})\times \bld\HH(\bar{K})\times \bld\HH(\bar{K})
\times \HH(\bar{K})\rightarrow\; 
H_3(K)\times E_3(K)\times V_3(K)\times W_3(K)$ 
are the following ones 
satisfying
\begin{alignat*}{2}
 \Pi_H^3 u(\vt) =&\; u(\vt)&&\quad \forall \vt\in \VT(K),\\
 \left(\grads \tr_H^\eg(\Pi_H^3 u),v
 \right)_\eg =&\; (\grads \tr_H^\eg(u),v)_\eg&&\quad \forall v\in 
 \grads\overset{\circ}{H_1}(\eg),
 \forall \eg\in \EG(K),\\
\left (\grads\tr_H^\fc(\Pi_H^3 u),v\right)_\fc =&\; 
\left (\grads\tr_H^\fc(u),v\right)_\fc&&\quad \forall 
 v\in\grads \overset{\circ}{H_2}(\fc), \forall \fc\in \FC(K),\\
\left (\grads \Pi_H^3 u,v\right)_K=&\; 
\left (\grads u,v\right)_K &&\quad \forall 
 v\in\grads \overset{\circ}{H_3}(K), 
\end{alignat*}
and
\begin{alignat*}{2}
\left(\tr_E^\eg(\Pi_E^3 u),v\right)_\eg =&\; \left(\tr_E^\eg(u),v\right)_\eg&&
\quad \forall v\in \oo{W_1}(\eg)\oplus\pol_0(K),
 \forall \eg\in \EG(K),\\
 \left(\curls\tr_E^\fc(\Pi_E^3 u),v\right)_\fc =&\;
  \left(\curls\tr_E^\fc(u),v\right)_\fc&&
 \quad \forall v\in \curls\overset{\circ}{E_2}(\fc),
 \forall \fc\in \FC(K),\\
  \left(\tr_E^\fc(\Pi_E^3 u),v\right)_\fc =&\; 
   \left(\tr_E^\fc(u),v\right)_\fc&&
 \quad \forall v\in \{v\in \overset{\circ}{E_2}(\fc):\;\curls v=0\},
 \forall \fc\in \FC(K),\\
 \left(\curls \Pi_E^3 u,v\right)_K =&\;
  \left(\curls u,v\right)_K&&
 \quad \forall v\in \curls\overset{\circ}{E_3}(K),\\
  \left( \Pi_E^3 u,v\right)_K =&\; 
   \left( u,v\right)_K&&
 \quad \forall v\in \{v\in \overset{\circ}{E_3}(K):\;\curls v=0\},
\end{alignat*}
and
\begin{alignat*}{2}
\left(\tr_V^\fc(\Pi_V^3 u),v\right)_\fc =&\;
  \left(\tr_V^\fc(u),v\right)_\fc&&
 \quad \forall v\in \oo{W_2}(\fc)\oplus\pol_0(\fc),\forall \fc\in \FC(K),\\
 \left(\divs \Pi_V^3 u,v\right)_K =&\;
  \left(\divs u,v\right)_K&&
 \quad \forall v\in \divs\overset{\circ}{V_3}(K),\\
  \left( \Pi_V^3 u,v\right)_K =&\; 
   \left( u,v\right)_K&&
 \quad \forall v\in \{v\in \overset{\circ}{V_3}(K):\;\divs v=0\},
\end{alignat*}
and
\begin{alignat*}{2}
 (\Pi_W^3 u,v)_K =&\; (u,v)_K&&\quad \forall v\in \oo{W_3}(K)\oplus\pol_0(K).
\end{alignat*}

\

\section*{Appendix B: proofs of sequence exactness}
\subsection*{Proof of Lemma \ref{lemma:2d-exact}}
\begin{proof}
We only provide a sketch of the proof since it  is very simple.
The exactness of these four  sequences 
can be proven using exactly the same argument by 
showing  that the differential operators map the previous 
function space into the next one, 
the $\mathrm{curl}$ operator is surjective, 
and  the following dimension count
holds
\[
 \dim H_2(K)-\dim E_2(K)+\dim W_2(K) = 1,
\]
where $H_2(K)$, $E_2(K)$, $W_2(K)$ are the related spaces
for the sequence. This completes the sketch of the proof.
\end{proof}

\subsection*{Proof of Lemma \ref{lemma:3d-exact}}
\begin{proof}
The exactness of the first two sequences, $\mathrm{S_{1,k}^{2d}}$ and 
$\mathrm{S_{2,k}^{2d}}$, is well-known; 
see \cite{ArnoldFalkWinther06}
for an elegant proof in arbitrary space dimensions which
includes the three-dimensional result as a special case, and uses the {\it Koszul complex} and the
{\it homotopy formula}.

Next, we prove that the last two sequences, 
$\mathrm{S_{5,k}^{2d}}$ and 
$\mathrm{S_{6,k}^{2d}}$, are exact. We omit the proofs for the 
third and fourth sequences since they are similar to those we are going to present now.

We prove the exactness by showing that the following three
identities hold
\begin{subequations}
\label{exact3d}
\begin{align}
\label{e3d1}
\bR= &\; \kerl_{\grads}H_3,\\
\label{e3d2}
\grads H_3= &\; \kerl_{\curls}E_3,\\
\label{e3d3}
\divs V_3= &\; W_3,
\end{align}
 and that the following dimension count identity holds
\begin{equation}
\label{e3d4}
\dim H_3 -\dim E_3+\dim V_3-\dim W_3 = 1.
\end{equation}
\end{subequations}
The other equality, namely,
$ \curls E_3=  \kerl_{\divs}V_3,$
is a direct consequence of the above results.
To see this, we have 
\begin{alignat*}{2}
\dim \curls E_3 =&\; \dim  E_3 - \dim \kerl_{\curls}E_3
&& \text{ }
\\
=&\; \dim  E_3 - \dim \grads H_3
&& \quad\text{ by \eqref{e3d2}}\\
=&\; \dim  E_3 - \dim H_3 +1
&& \quad\text{ by \eqref{e3d1}}\\
=&\; \dim  V_3 - \dim W_3
&& \quad\text{ by \eqref{e3d4}}\\
=&\; \dim  V_3 - \dim \divs V_3
&& \quad\text{ by \eqref{e3d3}}\\
=&\; \dim  \kerl_{\divs}V_3.
\end{alignat*}

We start with the verification for the sequence $\mathrm{S_{6,k}^{2d}}$. 
The first equality \eqref{e3d1} is trivially satisfied. 
We also have 
\[
\divs\left(\!\!\!\! \begin{tabular}{c}
$0$\\
$0$\\
$\pol_k(x,y)\oplus \pol_{k+1}(z)$\\
\end{tabular}\!\!\!\!\right)
=\pol_k(x,y)\oplus \pol_{k}(z),
\]
which implies the third equality \eqref{e3d3}.
The dimension equality \eqref{e3d4} is also easy to verify by the fact that
\begin{alignat*}{2}
\dim H_3 = &\;
\dim \pol_{k+1|k+1} \\
=&\;
(k+2)^2(k+3)/2\\
\dim E_3=&
\;
\dim \bld{N}_k(x,y)\cdot \dim \pol_{k+1}(z)
+\dim \pol_{k+1|k}\\
=&\;
3(k+1)(k+2)(k+3)/2
\\
\dim V_3 = &\;
\dim \bld{RT}_k(x,y)\cdot \dim \pol_{k}(z)
+
\dim \pol_{k|k+1}\\
=&\;
(k+1)^2(k+3)+(k+1)(k+2)^2/2\\
\dim W_3 = &\;
\dim \pol_{k|k}\\
=&\;
(k+1)^2(k+2)/2.
\end{alignat*}
Now, we are left to prove the identity \eqref{e3d2}.
Since $\grads H_3\subset \kerl_{\curls} E_3$,
we just need to show that $\kerl_{\curls} E_3\subset \grads H_3$.
To this end, let $p$ be a function in $\kerl_{\curls} E_3$. 
Since $E_3 \subset \bpol_{2k+2}
$, we have 
$p=\grads q$ for a scalar polynomial function $q\in\pol_{2k+3}$.
We have 
\[
\left(\!\!\!\begin{tabular}{c}
 $ \partial_x p$\\
   $\partial_y p$
\end{tabular}\!\!\!\right)
 \in \bld{N}_k(x,y)\otimes \pol_{k+1}(z),\;\;
 \text{ and }
 \partial_z p\in \pol_{k+1}(x,y)\otimes\pol_k(z)
\]
Now, we write $q$ in terms of a polynomial in $z$ with its coefficients being polynomials
in $x$ and $y$:
\[
 q = \sum_{\alpha=0}^{2k+3} f_\alpha(x,y)z^\alpha.
\]
We have 
\[
\partial_z q = \sum_{\alpha=1}^{2k+3} \alpha f_\alpha(x,y)z^{\alpha-1} \in 
\pol_{k+1}(x,y)\otimes\pol_k(z).
\]
This implies that 
$f_\alpha(x,y)\in \pol_{k+1}(x,y)$ for $1\le \alpha\le k+1,
$
and 
$f_\alpha(x,y)=0$  for $\alpha \ge k+2.
$
Hence,
\[
 q = f_0(x,y)+\sum_{\alpha=1}^{k+1} f_\alpha(x,y)z^\alpha.
\]
Since $f_\alpha(x,y)\in \pol_{k+1}(x,y)$ for $1\le \alpha\le k+1$ and
$H_3:=\pol_{k+1}(x,y)\otimes\pol_{k+1}(z)$, we have
\[
 \sum_{\alpha=1}^{k+1} f_\alpha(x,y)z^\alpha\in H_3.
\]
Then, we have
\[
\left(\!\!\!\begin{tabular}{c}
 $ \partial_x f_0(x,y)$\\
   $\partial_y f_0(x,y)$
   \end{tabular}\!\!\!\right)
 \in \bld{N}_k(x,y),
 \]
which implies $f_0(x,y)\in \pol_{k+1}(x,y)$, hence $q\in H_3$. This completes the proof
of the equality \eqref{e3d2}. Hence the sequence $\mathrm{S_{6,k}^{2d}}$ is exact.

Now, we prove equalities \eqref{exact3d} for the sequence 
$\mathrm{S_{5,k}^{2d}}$.
The first equality \eqref{e3d1} is trivially satisfied. The third equality is due to the fact that
\[
\divs\left(\!\!\!\! \begin{tabular}{c}
$\pol_k(x,y)\oplus \pol_{k}(z)$\\
$0$\\
$\widetilde\pol_k(x,y)\oplus \pol_{k+1}(z)$\\
\end{tabular}\!\!\!\!\right)
=\pol_k(x,y)\oplus \pol_{k}(z).
\]
The dimension equality \eqref{e3d4} is also easy to verify by the fact that
\begin{alignat*}{2}
\dim H_3 = &\;
k+3+\dim \pol_{k|k} 
\\
\dim E_3=&
\;
2k+3+\dim \bpol_{k|k}
\\
\dim V_3 = &\;
k+1+\dim \bpol_{k|k}\\
\dim W_3 = &\;
\dim \pol_{k|k}.
\end{alignat*}
Now, we are left to prove the identity \eqref{e3d2}.
Again, we prove that $\kerl_{\curls} E_3\subset \grads H_3$.
Since the spaces in $\mathrm{S_{5,k}^{2d}}$ is included in the 
related spaces in $\mathrm{S_{6,k}^{2d}}$, which is an exact sequence,
we have any function $p\in \kerl_{\curls} E_3$ is a gradient of a function 
$q\in \pol_{k+1|k+1}$. 
We have
\[
\left(\!\!\!\begin{tabular}{c}
 $ \partial_x p$\\
   $\partial_y p$
\end{tabular}\!\!\!\right)
 \in \left(\!\!\!\begin{tabular}{c}
 $ \pol_{k|k}$\\
   $\pol_{k|k}$
\end{tabular}\!\!\!\right)
\oplus\left(\!\!\!\begin{tabular}{c}
 $y$\\
$-x$
\end{tabular}\!\!\!\right)\widetilde{\pol}_k(x,y)
,\;\;
 \text{ and }
 \partial_z p\in \pol_{k|k}\oplus\widetilde\pol_{k+1}z^k.
\]

Now, we show that the function $q$ is actually a function in the space
\[
H_3 = \pol_{k|k}\oplus\widetilde\pol_{k+1}(x,y)\oplus \{z^{k+1}\}.
\]

Again, we express $q$ as a polynomial of the variable $z$ with coefficients polynomials of $x$ and $y$:
\[
 q=\sum_{\alpha=0}^{k+1}f_\alpha(x,y)z^\alpha,
\]
where $f_\alpha(x,y)\in\pol_{k+1}(x,y)$.
Using the fact that 
$ \partial_z q\in \pol_{k|k}\oplus\widetilde\pol_{k+1}z^k$, we immediately get 
$f_\alpha(x,y)\in\pol_k(x,y)$ for $1\le \alpha\le k$.
Moreover, since $\partial_x q\in \pol_{k|k}\oplus y\widetilde\pol_k(x,y)$, we have 
$\partial_x f_{k+1}(x,y) = 0$. Similarly, $\partial_y f_{k+1}(x,y) = 0$. This implies that 
$f_{k+1}(x,y)$ is a constant. Hence, $q\in H_3$ as desired. 
This completes the proof that $\mathrm{S_{5,k}^{2d}}$
is an exact sequence
and  completes the proof of Theorem \ref{lemma:3d-exact}.
\end{proof}


\end{document}